\documentclass[10pt]{article}

\usepackage{slashbox}
\usepackage[utf8]{inputenc}
\usepackage{a4wide}
\usepackage{bm}
\usepackage{graphicx}
\usepackage{amsmath}
\usepackage[parfill]{parskip}
\usepackage{xcolor}
\usepackage[english]{babel}
\usepackage{amsfonts}
\usepackage{amssymb}
\usepackage{amsbsy}
\usepackage{hyperref}
\usepackage{epstopdf}
\ifpdf
  \DeclareGraphicsExtensions{.eps,.pdf,.png,.jpg}
\else
  \DeclareGraphicsExtensions{.eps}
\fi

\newcommand{\dd}{\,\mathrm{d}}
\newcommand{\deltaUnstab}{{\delta_{||}}}
\newcommand{\deltaStab}{{\delta_\perp}}

\newtheorem{lemma}{Lemma}
\newtheorem{Corollary}{Corollary}

\newtheorem{theorem}[lemma]{Theorem}
\newtheorem{assumption}{Assumption}

\newcommand{\R}{\mathcal{R}}
\newcommand{\B}[1]{\mbox{\boldmath${#1}$\unboldmath}}
\newcommand{\D}[1]{{#1}'}
\newcommand{\DD}[1]{{#1}''}
\newcommand{\X}{\mathcal{X}}
\newcommand{\PP}{\mathcal{P}}
\newcommand{\PPplus}{\widehat{\PP}}
\newcommand{\QQ}{\mathcal{Q}}
\newcommand{\QQplus}{\widehat{\QQ}}
\newcommand{\blockdiag}[1]{\ensuremath \text{blockdiag}\,({#1})}

\DeclareMathOperator{\Span}{span}

\newcommand{\TheTitle}{Projected Shadowing-based Data Assimilation}

\title{\TheTitle}

\author{
Bart de Leeuw\thanks{Centrum Wiskunde \& Informatica, PO Box 94079, 1090 GB Amsterdam, Netherlands, 
	({b.m.de.leeuw@cwi.nl}, {s.dubinkina@cwi.nl}).} 
\and 
Svetlana Dubinkina\footnotemark[2]
\and 
Jason Frank\thanks{Utrecht University, Mathematical Insitute, P.O.~Box 80010, 3508 TA Utrecht, Netherlands, 
	({j.e.frank@uu.nl}).}
\and 
Andrew Steyer\thanks{Sandia National Laboratories, PO Box 5800 Albuquerque, NM 87185, 
	({asteyer@sandia.gov}).}
\and 
Xuemin Tu\thanks{University of Kansas, Dept. of Mathematics, 405 Snow Hall, Lawrence, KS 66405, 
	({xuemin@ku.edu}, {erikvv@ku.edu}).}
\and 
Erik Van Vleck\footnotemark[5]
}

\usepackage{amsopn}

\begin{document}
\maketitle

\begin{abstract}
In this article we develop algorithms for data assimilation based upon a computational time dependent stable/unstable splitting. Our particular method 
is based upon shadowing refinement and synchronization techniques and is motivated by work on Assimilation in the Unstable Subspace (AUS) \cite{CaGhTrUb08, TrDiTa10, PaCaTr13} and Pseudo-orbit
Data Assimilation (PDA) \cite{JuSm01,JuReRoSm08,DuSm14}. 
The algorithm utilizes time dependent projections onto the non-stable subspace determined by employing computational techniques
for Lyapunov exponents/vectors. The method is extended to parameter estimation without changing the problem dynamics and we address techniques for
adapting the method when (as is commonly the case) observations are not available in the full model state space.
We use a combination of analysis and numerical experiments (with the Lorenz 63 and Lorenz 96 models) to illustrate the efficacy of the techniques
and show that the results compare favorably with other variational techniques.
\end{abstract}

\begin{center}
\begin{textbf}
{Key words:}\end{textbf}
Data assimilation; tangent space decomposition; shadowing; synchronization; 
\end{center}

\begin{center}
\begin{textbf}
{AMS subject classifications:}\end{textbf}
62M20; 37C50; 34D06; 37M25
\end{center}

\section{Introduction\label{sec:introduction}}
Data assimilation methods combine orbits from a dynamical system model with measurement data to obtain an improved estimate for the state of a physical system.   
In this paper we develop a data assimilation method in the context of a discrete deterministic model
\begin{equation}\label{map}
	x_{n+1} = F_n (x_n), \quad x_n \in \R^d, \quad n=0,\dots,N-1,
\end{equation}
where $F_n: \R^d \to \R^d$. 
In many applications the model is defined by the time-discretization of an ordinary differential equation (ODE)
\begin{equation}\label{ode}
	\dot{x} = f(t, x), \quad x(t) \in \R^d,
\end{equation}
which in turn may be defined as the space-discretization of a partial differential equation (or system of PDEs).  

Let the sequence\footnote{In the sequel we will adopt the notation $\{ \X_n; n=0,\dots, N\}$ or simply as $\{ \X_n \}$ for a discrete orbit.  The latter notation is also occasionally employed to denote an infinite sequence.} 
$\{\X_0, \X_1, \dots, \X_N\}$ be a distinguished orbit of \eqref{map}, referred to as the \emph{true solution} of the model, and presumed to be unknown.
Suppose we are given a sequence of noisy observations $y_n$  related to $\X_n$ via 
\begin{equation}\label{obs}
	y_n = H \X_n + \xi_n, \qquad y_n \in \R^b, \quad n=0,\dots,N,
\end{equation}
where $H:\R^d \rightarrow \R^b$, $b\le d$, is the observation operator, and the noise variables $\xi_n$ are drawn from a normal distribution $\xi_n \sim \mathcal{N}(0,E)$ with zero mean and known observational error covariance matrix $E$.  

Data assimilation is the problem of finding an orbit (or pseudo-orbit, see \ref{sec:shadowing}) $\B{u}=\{ u_0, u_1, \dots, u_N\}$, $u_n \in \R^d$, of the model \eqref{map}, such that the differences $\| y_n - H u_n\|$, $n=0,\dots,N$ are small in an appropriately defined sense. This is done with the aim of minimizing the unknown error $\| u_n - \X_n \|$; see for example \cite{Ta97,LaStZy15}.  
For instance, well known four-dimensional variational data assimilation (4DVar) aims at finding the optimal initial condition $u_0$ of \ref{map} to minimize a cost function
\[
 C_\mathrm{var}(u_0; \{ y_n \}) = 
\sum_{n=1}^N (y_n-Hu_n)^T E^{-1}(y_n-Hu_n) + \lambda_n(u_n - F_n(u_{n-1})),
\]
where the $\lambda_n$ are Lagrange multipliers to ensure that the sequence $\{ u_n \}$ defines an orbit of \ref{map} (see e.g.\cite{Sa70,LeDe85,TaCo87,Ta97} and references therein). 
One drawback of variational data assimilation is that the number of local minima of the cost function increases dramatically with $N$ \cite{Be91,MiGhGa94,PiVaTa96}.
This places a practical limit on the length of the assimilation window---the time period over which observations may be assimilated. 

We propose a novel data assimilation method that overcomes this drawback: with the proposed method, increasing the length of the assimilation window may in fact lead to a better estimation. 
Instead of minimizing a cost function, we search for a zero of the \emph{cost operator}
\begin{equation}\label{residual}
	G(\B{u}) = \begin{pmatrix} G_0(\B{u}) \\ G_1(\B{u}) \\ \vdots \\ G_{N-1}(\B{u}) \end{pmatrix},
	\qquad G_n(\B{u}) = u_{n+1} - F_n(u_n), \quad n=0,\dots,N-1,
\end{equation}
using a contractive iteration started from (a proxy of) complete, noisy observations. 
This approach is motivated by research on numerical shadowing methods. We stress that, as is the case with 4DVar, our approach attempts to find an exact orbit of \ref{map} consistent with the observations.  However, instead of solving directly for the initial condition, we solve for the whole orbit at once.

As stated, our approach assumes the availability of (noisy) observations of the complete state vectors $\X_n$.  In other words, we assume that the observation operator $H$ is the identity matrix on $\R^d$. When only partial observations are available, it is necessary to generate a proxy for complete observations.  This can be done by some other cheap but inaccurate data assimilation method. For instance, in \ref{Sec:L63} we demonstrate this idea using direct insertion of noisy partial observations into the iteration \ref{map}.  

Recent attempts \cite{CaGhTrUb08, TrDiTa10, ToHu13, PaCaTr13, LaSaShSt14, SaSt15} to improve speed and reliability of data assimilation specifically address the partitioning of the tangent space into stable, neutral and unstable subspaces corresponding to Lyapunov vectors associated with negative, zero and positive Lyapunov exponents, respectively (see \ref{sec:QR}). 
In particular, Trevisan, d'Isidoro \& Talagrand propose a modification of 4DVar, so-called 4DVar-AUS, in which corrections are applied only in the unstable and neutral subspaces \cite{TrDiTa10,PaCaTr13}.
On the other hand, the stable subspace may also be exploited. Research by Pecora \& Carroll \cite{PeCa90} indicates that when partial observations are sufficient to constrain the unstable subspace, an orbit of the chaotic Lorenz 63 system can be made to converge exponentially in time to a different, driving orbit.
Their work has triggered a substantial body of research on the idea of synchronization of chaos 
(see review articles by Pecora et al.~\cite{PeCaJoMaHe97} and Boccaletti et al.~\cite{BoKuOsVaZh02}).

Motivated by the above, in this paper we propose a new method for data assimilation that utilizes distinct treatments of the dynamics in the stable and non-stable directions. 
We find a numerical orbit compatible with observations by using Newton's method with updates projected on the non-stable subspace to emphasize the need to stay close to 
current observations in non-stable directions. 
In the stable subspace, we ensure that the trajectory is determined by past observations using a forward integration to synchronize the stable components.
Although our focus here is on splitting into non-stable and stable components and then applying shadowing refinement and synchronization techniques, respectively, the splitting framework allows for other possibilities. In particular, if the non-stable subspace is relatively low dimensional this makes applications of techniques such as particle filters appealing. In addition, 4DVar or Kalman filter techniques may be applied to the stable system with the advantage that these techniques are being applied to a system with contractive dynamics. This also allows the split system to be put in a Bayesian data assimilation context.

In the next section we provide relevant background results. In \ref{sec:fullnewton} we describe the sense in which Newton's method is an effective data assimilation algorithm. While effective, the full Newton's iteration can be made more efficient by restricting the updates to just the non-stable tangent directions, as described in \ref{sec:projNewt}. The updates can then be synchronized in the stable directions as shown by the analysis in \ref{sec:analysis}. We provide details of our implementation in \ref{sec:implement}. Finally, in \ref{sec:num} we provide numerical results for the Lorenz 63 model and compare the method to 4DVar for the Lorenz 96 model.  We draw conclusions in \ref{sec:conclusions}.

\section{Background\label{sec:background}}

In \ref{sec:shadowing} we review concepts from numerical shadowing, in \ref{sec:QR} we describe the computation of tangent space splitting used in this paper, and in \ref{sec:synch} we review synchronization of chaos.

\subsection{Numerical shadowing\label{sec:shadowing}}

An $\varepsilon$-pseudo-orbit is a sequence $\B{u} =\{ u_0, u_1, \dots, u_N \}$ satisfying $\|G_n(\B{u})\| < \varepsilon$, $n=0,\dots,N-1$.  
For instance, suppose $F\equiv F_n$ is the exact time-$\tau$ flow map of an autonomous ODE $\dot{x}=f(x)$. 
If the components of $\B{u}$ are the iterates of a numerical integrator with local truncation error bounded by $\varepsilon$, then these define an $\varepsilon$-pseudo-orbit of $F$.  
The shadowing lemma (e.g.~Theorem 18.1.2 of \cite{KaHa95}) states that in a neighborhood of a hyperbolic set for $F$, for every $\delta>0$ there exists $\varepsilon>0$ such 
that every $\varepsilon$-pseudo-orbit is $\delta$-shadowed by an orbit of $F$, i.e.~there exists an orbit $\{ x_n \}$ satisfying $x_{n+1}=F(x_n)$ such that $\| u_n - x_n \|<\delta$ for all $n=0,\dots,N$.  
Rigorous bounds on the global error of numerical integrations with respect to a shadowing orbit can be proved by applying the Newton-Kantorovich theorem to Newton's iteration for $G(\B{x}) = 0$ with 
starting data given by the numerical iterates $\B{u}$ on a time interval that is long relative to the characteristic Lyapunov time
\cite{Be87, HaYoGr87, HaYoGr88, ChLiPa89, ChPa92, VV01}.
Shadowing is an important analysis technique for obtaining global error bounds on the numerical approximation to the solution of differential equations exhibiting chaos. 
We can view data assimilation in the same vein by interpreting the data as some approximation to the model solution and set it as our goal to find a particular model solution that shadows the data.

With respect to shadowing, the \emph{inverse problem} is to determine an optimal initial condition $u_0$ for a numerical integration, 
such that the numerical iterates 
$\B{u}$ $\delta$-shadow a desired orbit of \eqref{ode}.    
Shadowing refinement (see, e.g., \cite{GHYS90}), employs the pseudo-orbit as an initial guess for $G(\B{u})=0$ and, as opposed to proving the existence of a nearby zero of $G$, iteratively
refines the pseudo-orbit to obtain an improved approximation of a true solution.  This is clearly akin to the data assimilation problem.\par

Shadowing theory has already motivated a practical data assimilation algorithm known as pseudo-orbit data assimilation (PDA); see for instance \cite{JuSm01,JuReRoSm08,DuSm14} and references therein.
For the PDA approach a cost function 
\[
 C_\mathrm{PDA} = \frac{1}{2} \sum_{n=0}^{N-1} G_n^T G_n
\]
is minimized and the minimization is also initialized from observations. Obviously, the (nonunique) global minimum of $C_\mathrm{PDA}$ is zero and this value is reached if and only if $G(\B{u})=0$, that is, 
if $\B{u}$ is any model trajectory. The approach in \cite{DuSm14} approximatly minimizes $C_\mathrm{PDA}$ by taking a fixed number of gradient descent steps starting from observations. 
This typically yields not an orbit but a (discrete) pseudo-orbit, i.e.~the minimizing sequence satisfies $\| u_n - \X_n \| < \varepsilon$, for all $n=0,\dots,N$, and some constant $\varepsilon$. The distance between the pseudo-orbit and the manifold of trajectories is then smaller than the distance between observations 
and the manifold of trajectories. The mid-point of this pseudo-orbit is then used as the initial condition for a trajectory that should be consistent with model and data. 
PDA has been applied in operational weather models \cite{JuReRoSm08}, parameter estimation \cite{SmCuDuJu10} and as a method for finding reference trajectories for ensemble forecasting \cite{DuSm14}.

\subsection{Tangent subspace decomposition\label{sec:QR}}

In this section we review the decomposition of the tangent space into stable, neutral and strongly unstable subspaces.  This decomposition is central to the method described in this paper.  Let $\{ x_n; n=0,\dots,N \}$ denote an orbit of \ref{map}.  The fundamental matrix equation associated with $\{ x_n \}$ is a matrix valued difference equation
\begin{equation}\label{fundmateq}
	X_{n+1} = \D{F}_n (x_n) X_n, \quad n=0,\dots,N-1,
\end{equation}
where $X_n \in \R^{d\times d}$.
The iterates of \ref{fundmateq} become increasingly ill-conditioned as the columns align with the dominant growth direction.  To stably estimate $X_n$, one may introduce a time-discrete QR factorization. Let $X_0=Q_0R_0$,
and write 
\begin{equation}\label{Lyamgs}
Q_{n+1}R_{n+1} = \D{F}_n(x_n)Q_n \quad \mbox{for} \quad n=0,...,N-1,
\end{equation}
where $\D{F}_n(x_n)Q_n$ is a matrix product of known quantities, and $Q_{n+1}R_{n+1}$
is the $QR$ factorization found using the modified Gram-Schmidt process.
Then $X_1 = \D{F}_0(x_0)Q_0R_0 = Q_1R_1R_0$, $X_2 = \D{F}_1(x_1)X_1 = \D{F}_1(x_1)Q_1R_1R_0 = Q_2R_2R_1R_0$, etc.
Note that this procedure is well defined for $Q_n\in\R^{d\times p}$ for $p\leq d$ provided
$\D{F}_n(x_n)Q_n$ is full rank for all $n$. 
The Gram-Schmidt process yields the unique upper triangular $R_n\in \R^{p\times p}$ with positive diagonal elements and, importantly, preserves the ordering of the columns of the $Q_n$.

The (local) $p$ ($1\leq p\leq d$) largest Lyapunov exponents of the orbit $\{x_n \}$ are extracted from the time average of the logarithm of the diagonal of $R_n$ \cite{DiVV15}:
\begin{equation*}
	\lambda_i = \limsup_{N\rightarrow\infty} \frac{1}{N} \sum_{n=1}^N \ln R_n^{(i,i)}, \quad i=1,\dots,p.
\end{equation*}
The method of construction ensures $\lambda_1 \ge \lambda_2 \ge \cdots \ge \lambda_p$.  
Associated with $\lambda_i$ is a Lyapunov vector $V_n^{(i)}$. The columns of $Q_n$ generally (for most initial conditions)
form an orthonormal basis for the Lyapunov vectors at time $n$ (see \cite{DE06,DEVV10,DEVV11}). The iteration \ref{fundmateq} is a generalized power iteration. For each $\ell=1,\dots,p$, one finds in the limit 
$n\rightarrow \infty$ that $\Span \{ V_n^{(1)}, \dots,V_n^{(\ell)} \} = \Span \{ Q_n^{(1)}, \dots, Q_n^{(\ell)} \}$, where $Q_n^{(i)}$ denotes the $i$th column of $Q_n$ \cite{DiVV15, DiVV05, Ad95, BeGaGiSt80}. 
Positive (negative) $\lambda_i$ correspond to tangent directions $Q^{(i)}_n$ in which perturbations grow (decay) exponentially.   Consequently, if $\lambda_{p} \ge 0 > \lambda_{p+1}$, 
then the matrix $Q_n^u = ( Q_n^{(1)}, \dots  Q_n^{(p)})$ provides an orthonormal basis for the non-stable tangent space at $X_n$. In practice we may obtain $Q_n^u$ using a thin QR-factorization,
\begin{equation}\label{Lyamgsthin}
Q_{n+1}^uR_{n+1}^u = \D{F}_n(x_n)Q_n^u \quad \mbox{for} \quad n=0,...,N-1,
\end{equation}
where $Q_0^u\in\R^{d\times p}$ and $R_n^u\in \R^{p\times p}$ (please note $R_n^u$ is the upper left $p\times p$ block of $R_n$). Note that this procedure is well defined for $Q_n\in\R^{d\times p}$ for $p\leq d$ provided
$\D{F}_n(x_n)Q_n$ is full rank for all $n$. 
The Gram-Schmidt process yields the unique upper triangular $R_n^u\in \R^{p\times p}$ with positive diagonal elements and, importantly, preserves the ordering of the columns of the $Q_n^u$.
We note here that by approximating the non-stable subspace 
we obtain information (see \cite{ToHu13}) that may be used to
analyze the error in data assimilation schemes, namely in terms of the degree to which
observations constrain the uncertainty within the non-stable subspace.  
We remark that the dimension of the unstable subspace may be much less then the total dimension. In Carrassi et al. \cite{CaTrDeTaUb08} it is shown that the AUS-framework gives good results for a quasi-geostrophic model described in \cite{RoBa96}. This model is of dimension 14784, while the unstable subspace has a dimension of 24 \cite{SnHa03}.

We will use the computed factors $Q_n\in\R^{d\times p}$ to construct projection operators onto the non-stable tangent space. 
The $Q_n$ are quantities that can be computed robustly with good forward error analysis properties (under reasonable assumptions closely
related to the continuity of Lyapunov exponents with respect to perturbations). In particular, the results in \cite{DVV06,DVV08,VV09,BVV12}
show that the $Q_n$ are continuous with respect to errors in $F'(x_n)$ and quantify the error in the $Q_n$
as a function of the separation in growth/decay rates. 
This is characterized by the integral separation or integral separation structure (see also \cite{Ad95}) which is
closely related to the continuity of Lyapunov exponents with respect to
perturbations of $F'(x_n)$. 
In our context this ensures the time dependent projection operators $P_n = Q_n^u Q_n^u T$ are robust. 

\subsection{Synchronization\label{sec:synch}}

Pecora \& Carroll \cite{PeCa90} demonstrated that an orbit of a chaotic dynamical system (the observer) can sometimes be made to synchronize with a second orbit (the driver) of that system, given partial observations of the driver signal. There is a sizeable body of literature on synchronization of chaos, particularly in the field of systems and control \cite{PeCa90,HaOlTi11,ToHu13}.

\begin{figure}[!htb]
\begin{center}\includegraphics[width=0.8\textwidth]{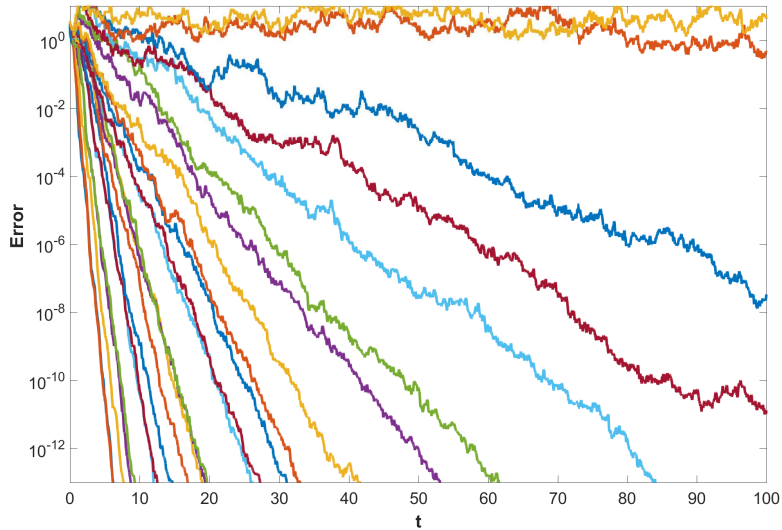}
\end{center}
\caption{Synchronization in the L96 model \ref{L96} \cite{Lo96} with 13 positive Lyapunov exponents in 36 dimensions. We plot the $\ell^\infty$-norm of the difference between the true solution and the synchronization approximation as a function of time. Forcing is done with projections of the true solution onto the non-stable space and the different graphs are for $k=12,13,...,35$, colored with the Matlab default color order starting from 35. It can be observed that after a transient time and for sufficiently large $k$ (i.e. $k\geq 14$), convergence to the true solution is exponential.}
\label{artikel_sync}
\end{figure}

For our purposes, the following coupled driver-response process is appropriate:
\begin{subequations}
\begin{align}
	x_{n+1} &= F_n(x_n), \label{driver} \\
	z_{n+1} &= P_n x_{n+1} + (I-P_n) F_n(z_n), \label{receiver}
\end{align}
\end{subequations}
where the $P_n\in \R^{d\times d}$ are a sequence of appropriately chosen projection matrices. The manifold $\mathcal{S} = \{ (x,z) \in\R^d \times \R^d : x=z\}$ is invariant under these dynamics and is called the \emph{synchronization manifold}.  When $\mathcal{S}$ attracts a neighborhood of itself, then for $z_0$ within the basin of attraction, $z_n$ synchronizes with $x_n$.  Defining $w_n = z_n - x_n$, $n=0,1,\dots,$ the transverse dynamics with respect to $\mathcal{S}$
is given by
\begin{align*}
w_{n+1} &= P_n x_{n+1} + (I-P_n) F_n(z_n) - F_n(x_n) \notag \\
&= P_n F_n(x_n) + (I-P_n) F_n(z_n) - F_n(x_n) \notag \\
&= (I-P_n) \left[ F_n(x_n+w_n) - F_n(x_n) \right] \notag \\
&= (I-P_n) \D{F}_n (x_n) w_n + r_n(w_n), \notag
\end{align*}
where $r_n(w)$ is assumed to be of higher order in $w$. The projectors $P_n$ need to be chosen to ensure asymptotic stability of the origin under the transverse dynamics.  From the stability theory of Lyapunov, it is known that if the sequence $\|w_n\|$ converges exponentially to zero for generic initial conditions,  then the Lyapunov exponents of the transverse dynamics must necessarily all be negative.  Such a necessary condition is argued by Pecora \& Carroll in \cite{PeCa91}.  On the other hand, negativity of the Lyapunov exponents is also sufficient for convergence in a neighborhood of the origin, if $\D{F}_n$ is regular and $r_n$ is at least second order in $w$.
In our application to data assimilation we will choose $P_n$ to project (in an approximate sense) onto the locally non-stable tangent space $Q_n^u$. \ref{artikel_sync} illustrates synchronization of the Lorenz 96 model (see \cite{Lo96} and \ref{Sec:L96}) using the driver-response system 
\eqref{driver}--\eqref{receiver} with projection $P_n=Q_n^u Q_n^u T$ for increasing dimension of the projection space $p$. In particular we observe exponential convergence only when $p$ is greater than or equal to the dimension of the nonstable space, with exponential rate of convergence increasing with $p$.

\section{Data assimilation via Newton's method\label{sec:fullnewton}}
In this section we discuss the use of Newton's method for data assimilation, a context in which it was first applied in \cite{BrPa01}.
An important property of Newton's method is its local nature: when the initial guess is sufficiently close to a zero, the iterates converge to that zero at a quadratic rate. 
This statement is made formal in the Kantorovich Theorem \cite{KaAk59}. 

Consequently, by analogy to the shadowing approach to global error estimation, we may construct a simple scheme for data assimilation by applying Newton's iterations to solve
\[
 G(\B{u}) = 0,
\]
where $G$ is defined in \eqref{residual} and starting data is provided by the noisy observations $\{ y_n: n=0,\dots,N\}$ with observation operator the identity $\{ H \X_n = \X_n : n=0,\dots,N\}$. 
(an assumption that can be relaxed, see \ref{sec:num}).

In the $k$th Newton's iteration we seek an update $\B{\delta}^{(k)}$ approximately solving
\begin{equation}\label{Newton1}
	G (\B{u}^{(k)} + \B{\delta}^{(k)}) = 0.
\end{equation}
We then update using $\B{u}^{(k+1)} = \B{u}^{(k)} + \B{\delta}^{(k)}$.
The solution to \eqref{Newton1} is approximated by iterating
\begin{equation}\label{fullNeweq}
\D{G}(\B{u}^{(k)})\B{\delta}^{(k)} = -G(\B{u}^{(k)}), \quad \B{u}^{(k+1)} := \B{u}^{(k)} + \B{\delta}^{(k)}
\end{equation}
to convergence. We remark that the function $G(\B{u})$ has a zero for every orbit of the model \eqref{map}. The function $G:\R^{dN}\to \R^{d(N-1)}$ has a $d(N-1)\times dN$ Jacobian with block structure
\[
	\D{G}(\B{u}) = \begin{bmatrix} 
	- \D{F}_0(u_0) & I \\
	& -\D{F}_1(u_1) & I \\
	& & \ddots & \ddots \\
	& & & -\D{F}_{N-1}(u_{N-1}) & I \end{bmatrix}.
\]
We solve each Newton's step using the right pseudoinverse of $\D{G}$, i.e.~$\D{G}^\dag = \D{G}^T (\D{G} \, \D{G}^T)^{-1}$, where the linear system involving the block tridiagonal matrix $\D{G}\, \D{G}^T$ is solved using a block tridiagonal solver. To distinguish this method from the projected method to be described in \ref{sec:projNewt}, we shall refer to it as the \emph{full Newton's method}.

The fact that Newton's method is a \emph{local} root-finding method proves useful. Initializing it with observations, we can expect to find a trajectory close to observations, provided the initial observational error is not too large \cite{BrPa01, Br10}.

\begin{figure}[!htb]
\begin{center}
\includegraphics [width=0.8\textwidth]{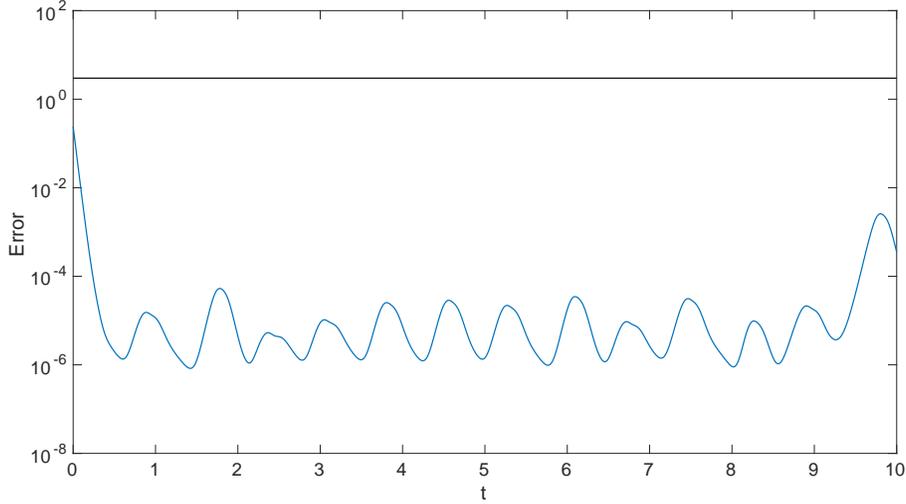}
\end{center}
\caption{Application of Newton's method to the L63 model \ref{L63}. The mean observational error is in black, the error over time of the estimation is in blue.}
\label{figfullNewton}
\end{figure}

The convergence of this approach with Newton's method can be demonstrated with a numerical example. 
Using the Lorenz 63 model \cite{Lo63} (see also \ref{Sec:L63}) the true trajectory $\{ X_n \}$ was integrated and perturbed\footnote{For this simple demonstration, we compute just a single realization of the noise process.  Later in \ref{Sec:L63} we include results for an ensemble.} as specified in \eqref{obs} with noise covariance 
$E = I$. The perturbed data was used as a starting guess for Newton's method. As shown in \ref{figfullNewton}, convergence to a model trajectory can be observed, 
and the mean square error (MSE) defined as
\begin{equation}\label{MSE}
 \mathrm{MSE} = 
\frac{1}{N}\sum_{n=1}^N (u_n-\X_n)^T (u_n-\X_n)
\end{equation}
is equal to 8.1921e-04. Next we examine errors with respect to the observation operator
\begin{equation}\label{anerror}
 C(\{x_n\}) = 
\frac{1}{N}\sum_{n=1}^N (y_n-Hx_n)^T (y_n-Hx_n).
\end{equation}
The mean squared error in the observations (mean noise variance) is given by $C(\left\lbrace \X_n \right\rbrace)=2.9853$.  By comparison, the mean observation discrepancy of the Newton's solution is equal to
$C(\B{u})=2.9844$.
That is, even though the trajectory found by Newton's method is not identical to the true trajectory, it is in fact a model orbit \emph{closer} 
to the observations. This demonstrates that the method works well, even when observational noise prevents determining a unique viable trajectory. 

\begin{table}
\caption{Application of Newton's method to L63. $C(\cdot)$ and MSE are defined in equations \eqref{anerror} and \eqref{MSE} respectively.}
\centering
\begin{tabular}{|l|c|}\hline
\textbf{Property}&\textbf{Value}\\\hline
Observation error $C(\left\lbrace \X_n \right\rbrace)$&2.9853\\
Observational discrepancy of estimated trajectory $C(\B{u})$&2.9844\\
MSE between estimation and the truth &8.1921e-04\\\hline
\end{tabular}
\end{table}
{\bf Remark.} Our data assimilation method can be applied to parameter estimation as well. A standard approach to dealing with parameter estimation is to treat parameters as dependent variables with trivial dynamics.  This approach adds neutral directions to the tangent space, which can hamper convergence of shadowing.  Instead,  
Newton's method can be extended to simultaneously estimate state space variables and parameters. 
Consider $G(\B{u})$ defined in (4) and replace $G_n(\B{u})$ with $G_n(\B{u};\bm{\alpha}) = u_{n+1} - F_n(u_n;\bm{\alpha})$ where $\bm{\alpha} = (\alpha_1, \alpha_2, ..., \alpha_q)$. Linearization with respect to
$\alpha_j$ takes the form
\[
\frac{\partial F_n(u_n;\bm{\alpha})}{\partial \alpha_j} d\alpha_j,\,\,\, n=0,...,N-1.
\]

In the presence of uncertain parameters, the linearization of $G$ is modified. In particular,
$\D{G}(\B{u})$ becomes $\D{G}(\B{u};\bm{\alpha})$. In the case of
$q$ parameters we have 
\[
\D{G}(\B{u};\bm{\alpha}) = [\D{G}_{\bm{u}}|\D{G}_{\bm{\alpha}}],
\]
where $G_{\bm{\alpha}}$ is composed of $q$ column vectors. 
Note that when forming the pseudoinverse of $\D{G}$, $\D{G}^\dag = \D{G}^T (\D{G} \D{G}^T)^{-1}$ we have
\[
\D{G} \D{G}^T = \D{G}_{\bm{u}} \D{G}^T_{\bm{u}} + \D{G}_{\bm{\alpha}} \D{G}^T_{\bm{\alpha}}
\]
so that $\D{G}_{\bm{\alpha}} \D{G}^T_{\bm{\alpha}}$ is a rank $q$ perturbation of
the block tridiagonal matrix $\D{G}_{\bm{u}} \D{G}^T_{\bm{u}}$. This allows for the use of Sherman-Morrision-Woodbury
formulas and a solver for $\D{G}_{\bm{u}} \D{G}^T_{\bm{u}}$ to solve linear systems with matrix $\D{G} \D{G}^T$.
In the case of time dependent parameters $\bm{\alpha}_n$, the diagonal blocks of $\D{G}_{\bm{u}} \D{G}^T_{\bm{u}}$ are modified and the overall
block tridiagonal structure in maintained in $\D{G} \D{G}^T$.  In \ref{sec:num_param_est} we illustrate this approach with numerical experiments.

\section{Tangent space splitting of Newton's method\label{sec:projNewt}}
In the previous section we demonstrated that Newton's method applied to the residual \ref{residual} may converge from noisy observations to a model trajectory. On the other hand the computational and memory costs of the full Newton's method may be high.  We will see that when the number $p$ of nonnegative Lyapunov exponents is moderate, substantial savings may be realized by computing Newton's updates only in the non-stable directions.

We start by decomposing the relation \eqref{Newton1} into the equivalent system
\begin{align*}
	\PP G (\B{u}^{(k)} + \PPplus \B{\delta}^{(k)} + (I-\PPplus) \B{\delta}^{(k)}) &= 0, \\
	(I-\PP) G (\B{u}^{(k)} + \PPplus \B{\delta}^{(k)} + (I-\PPplus) \B{\delta}^{(k)}) &= 0.
\end{align*}
Here, $\PP$ and $\PPplus$ are block diagonal projection matrices $\PP = \blockdiag{P_1,\cdots,P_N}$ and $\PPplus = \blockdiag{P_0,\dots,P_N}$,
where $P_0, P_1, ..., P_N\in\R^{d\times d}$ are projection matrices onto the non-stable subspace at time levels $n=0,1,\dots,N$, respectively. 

We propose to modify the Newton's iteration as follows. Instead of computing the update $\B\delta^{(k)}$ by simultaneously solving the above system, we split the iterate into updates in the range and complement of $\PPplus$. We also allow the projection operators $\PP$ and $\PPplus$ to be updated in each iteration.  In the $k$th iteration, we first approximate the update in the range of $\PPplus^{(k)}$, neglecting the term $(I-\PPplus^{(k)}) \B{\delta}^{(k)}$ in the first equation above and solving
\begin{equation}\label{unstcor}
	\PP^{(k)} G (\B{u}^{(k)} + \PPplus^{(k)} \B{\delta}^{(k)}) = 0
\end{equation}
for $\B{\deltaUnstab}^{(k)} = \PPplus^{(k)} \B{\delta}^{(k)}$.  Next we approximate the update in the complement of $\PPplus^{(k)}$ by solving
\begin{equation}\label{stabcor1}
	(I-\PP^{(k)}) G (\B{u}^{(k)} + \PPplus^{(k)} \B{\delta}^{(k)} + (I-\PPplus^{(k)}) \B{\delta}^{(k)}) = 0 
\end{equation}
for $\B{\deltaStab}^{(k)} = (I-\PPplus^{(k)})\B{\delta}^{(k)}$.  Then the update is computed as $\B{u}^{(k+1)} = \B{u}^{(k)} + \B{\deltaUnstab}^{(k)} + \B{\deltaStab}^{(k)}$.
Expressions \ref{unstcor} and \ref{stabcor1} are solved approximately for the components $\B{\deltaUnstab}^{(k)}$ and $\B{\deltaStab}^{(k)}$ as described below.

\subsection{Computation of projection matrices}
The basis $Q_n^u$, $n=0,\dots,N$, for the non-stable tangent space along the true trajectory $\{\X_n\}$ is unknown.  Instead, we approximate the $Q_n^u$ along the most recent approximate trajectory $\{ u_n^{(k)}\}$.  In each iteration we update the projection matrices $P_n^{(k)}$ that project 
onto the non-stable tangent space.  In the $k$th iteration We choose $P_n^{(k)} = Q_n^{u (k)} (Q_n^{u (k)})^T$, where $Q_n^{u (k)} \in \R^{d\times p}$ is a columnwise orthonormal matrix defined via the iteration \ref{Lyamgs} linearized along the most recently updated pseudo-orbit $\B{u}^{(k)}$.  That is, we take $x_n = u_n^{(k)}$, $n=0,\dots,N$, in \ref{Lyamgs}. For the first iteration we use the observations: $u_n^{(0)} = y_n$, $n=0,\dots,N$.

The dimension $p$ of the orthonormal basis $Q_n^{u (k)}$ should be equal to or greater than the number of non-negative Lyapunov exponents.    In practice we take $p$ to be a few more than the number of non-negative Lyapunov exponents to enhance the convergence rate of the synchronization step below \cite{DiRuVV97}.  

\subsection{Newton's step on the unstable space}
Linearization of \eqref{unstcor} yields a projected linear system for the update $\B{\deltaUnstab}^{(k)} = \PPplus^{(k)}\B{\delta}^{(k)}$:
\begin{equation}\label{projNewt}
	\PP^{(k)} \D{G}(\B{u}^{(k)}) \PPplus^{(k)} \B{\delta}^{(k)} = - \PP^{(k)} G(\B{u}^{(k)}). 
\end{equation}
Supressing the iteration index $k$ for the moment, define block matrices $\QQ = \blockdiag{Q_1^u,\dots,Q_N^u}$ and $\QQplus=\blockdiag{Q_0^u,\dots,Q_N^u}$, and note the relations $\QQ\QQ^T=\PP$, $\QQ^T\QQ = I$ with analogous expressions for $\QQplus$.   
Let $\B{\mu} = \QQplus^T \B{\delta} = \QQplus^T\PPplus\B{\delta}$, $\D{\widetilde{G}} = \QQ^T \D{G}(\B{u}) \QQplus$ and $\B{b} = \QQ^T G(\B{u})$.  Then the linear system for the update $\B{\mu}$ may be written as
\begin{equation} \label{reduced}
	\D{\widetilde{G}} \B{\mu} = -\B{b},
\end{equation}
where the matrix $\D{\widetilde{G}}$ has the block structure
\[
	\widetilde{G}' = \begin{bmatrix}  
	-R_0^u & I & & & \\
	& -R_1^u & I & & \\
	& & \ddots & \ddots & \\
	& & & -R_{N-1}^u & I 
	\end{bmatrix},
\]
and consequently, $\D{\widetilde{G}} \in \R^{Np \times (N+1)p}$. We solve \eqref{reduced} using the right pseudoinverse $\D{\widetilde{G}}^\dag = \D{\widetilde{G}}^T (\D{\widetilde{G}} \D{\widetilde{G}}^T)^{-1}$ 
and define the intermediate update 
\begin{equation}\label{update1}
	\B{\bar{u}}^{(k)} = \B{u}^{(k)} + \PPplus^{(k)} \B{\delta}^{(k)}
    = \B{u}^{(k)} + \QQplus^{(k)} \B{\mu}^{(k)}.
\end{equation}

\subsection{Synchronization step in the stable space}
We next turn to the treatment of \ref{stabcor1}.
Inserting the definition \ref{update1} into \ref{stabcor1} yields the relation
\begin{equation}\label{stabcor}
 	(I-\PP^{(k)}) G (\B{\bar{u}}^{(k)} + (I-\PPplus^{(k)}) \B{\delta}^{(k)}) = 0,
\end{equation}
whose solution for $\B{\deltaStab}^{(k)} = (I-\PPplus^{(k)})\B{\delta}^{(k)}$ we wish to approximate.  
Again dropping the iteration index $k$ for the moment, we expand \eqref{stabcor} component-wise over the time index $n$:
\begin{align*}
0&=\left[(I-\PP)G \left(\B{\bar{u}}  + (I - \PPplus)\B{\delta}\right)\right]_n\\
&=(I-P_{n+1})\left[\bar{u}_{n+1} + (I - P_{n+1})\delta_{n+1} - F_{n}(\bar{u}_{n} + (I - P_{n})\delta_{n})\right], \quad n=0,\dots,N-1.
\end{align*}
The second equation is rewritten in the form
\begin{equation*}
(I-P_{n+1})\delta_{n+1} = (I-P_{n+1})\left(F_{n}(\bar{u}_{n}+(I-P_n)\delta_{n})-\bar{u}_{n+1}\right), \quad n=0,\dots,N-1.
\end{equation*}
Adding $\bar{u}_{n+1}$ to both sides of this equation we get 
\[
	\bar{u}_{n+1} + (I-P_{n+1})\delta_{n+1} = P_{n+1}\bar{u}_{n+1} + (I-P_{n+1}) F_{n}(\bar{u}_{n}+(I-P_n)\delta_{n}), \quad n=0,\dots,N-1,
\]
or, defining $u_n^{(k+1)} = \bar{u}_n^{(k)} + (I-P_n) \delta_n^{(k)}$,
\begin{equation}\label{synceq}
	u_{n+1}^{(k+1)} = P_{n+1} \bar{u}_{n+1}^{(k)} + (I-P_{n+1}) F_n(u_n^{(k+1)}),\quad n=0,\dots,N-1.
\end{equation}
The form of this iteration is identical to that of the receiver equation \eqref{receiver} in the synchronization process. 
In other words, given the update $\B{\bar{u}}^{(k)}$, corrected in the non-stable subspace \ref{update1}, the correction to the stable subspace can be implemented through a forward synchronization integration \ref{synceq}.  In the Appendix, we prove that under suitable assumptions, if after $k$ iterations, the error in $\B{\bar{u}}^{(k)}$ exists entirely in the stable tangent space in the sense that  $P_n \bar{u}_n^{(k)} = \bar{u}_n^{(k)}$ and $\| P_n \mathcal{X}_n - \bar{u}_n^{(k)} \| < \epsilon$, then the forward integration \ref{synceq} converges exponentially to $\mathcal{X}_n$ as $n\to\infty$.

To summarize, the complete iteration step consists of:
\begin{enumerate}
\item Compute the approximate basis $Q_n^u$, $n=0,\dots,N$, for the tangent bundle along the pseudo-trajectory $\{u_n^{(k)}\}$.
\item Solve the linear system \ref{projNewt} for the update $\B{\deltaStab}^{(k)}$ in the non-stable subspace, and compute the intermediate update $\B{\bar{u}}^{(k)}$ from \ref{update1}.
\item Synchronize in the stable subspace using the forward iteration \ref{synceq} to obtain $\B{u}^{(k+1)}$.
\end{enumerate}

The Newton's step in the unstable subspace is based upon residual ($r_n := u_{n+1} - F_n(u_n)$) correction with both the residual and
the correction projected into the unstable subspace. If $P_nr_n = r_n$ for all $n$, i.e, the residual is wholly within the unstable subspace,
then the synchronization step in the stable subspace is trivial with $(I-P_n)\delta_n \equiv 0$
for all $n$. Thus, provided the Newton's iteration converges, all residual correction occurs within the unstable subspace.
In the more general case in which the residual is contained, at least for some $n$, in both the stable and unstable subspaces,
then the initialization of the synchronization step makes possible a reduction of the residuals in the stable subspace. This
then generates an updated approximate trajectory to linearize about and obtain updated projections. In this case the Newton's
step in the unstable subspace may again decrease the residual with respect to these new projections. The process then continues
in the updated stable subspace and we continue until the desired tolerance is achieved or the method fails for lack of convergence
of the projected Newton's iteration. In general the projected Newton's iteration will converge provided the residuals are small enough
as compared to the strength of the hyperbolic structure (exponential dichotomies, etc.) in the projected system.

What we have observed is that better results are obtained by switching after each projected Newton's
iterate to the synchronization step as opposed to switching to the synchronization after the projected Newton's has converged to tolerance.
We attribute this to the variation in the projections that are produced. We note here that the basic splitting based upon projection
into unstable and stable parts allows for different techniques to be employed for each subsystem. It provides a representation
for the unstable subspace which we believe will prove useful in assessing the effectiveness and uncertainties in data assimilation techniques. 
We also emphasize that in contrast with traditional data assimilation techniques but similar to PDA, the only influence of the observations 
is via the initial guess for the projected Newton's/synchronization scheme. Thus, in a perfect model scenario convergence to a solution depends on the initial guess being within its basin of attraction.

In the next section we demonstrate the algorithm and compare it to 4DVar for a number of test problems. 

\section{Implementation\label{sec:implement}}

In this section we provide details of the algorithm we implement and discuss some possible variations. The algorithm is ``interval sequential'' in the sense that the shadowing refinement is applied over an entire subinterval. This has the effect of simultaneously incorporating all observations over this subinterval into a single refinement step. In order to transition between subintervals we impose a continuity constraint in the stable subspace. Also discussed in this section are methods for obtaining an initial approximation of the solution trajectory. This is needed in order to determine the initial projections on each subinterval.

When observations are not available at every time step, we can redefine $F$ to be the map corresponding to the composition of several time steps (examples are given in \ref{Sec:L96,Sec:compto4dvar}). 
When observations are not of the full model state, we can first apply a preprocessing step to the observations to infer an estimate of the full state at all observation times and 
then perform the main algorithm with the goal of substantial noise reduction.
For the PDA method, where the same issue arises, this completion has been done using a variational analysis \cite{JuReRoSm08} or by just inserting climatological means for missing observations \cite{Du09,DuSm14}. 
In \ref{Sec:L63}, we demonstrate an alternative preprocess motived by synchronization, whereby the observation data is directly inserted as a driving signal. The effectiveness of such an approach relies on the ability of the partial observational data to constrain the unstable tangent space. However, it is one of the main conclusions of this paper that such a requirement on the data must hold anyway, if data assimilation is to be effective.

In our implementation, we decompose the time interval $t\in[0,T]$ of integration into $M+1$ non-overlapping time windows, and the data assimilation method is applied sequentially on each of these. We identify times $\tau_m$, $m=0,\dots,M+1$, where  $\tau_0=0$, $\tau_{M+1}=T$, $\tau_1\in(0,T)$ is the length of the first time window, and $\tau_m = \tau_1 + (m-1)\Delta \tau$, $\Delta \tau = (T-\tau_1)/M$.  The $m$th time window is the interval $t\in[\tau_{m-1},\tau_m]$.
In each window, an initial condition $(\deltaStab)_0$ is needed for the synchronization step, and convergence of the stable directions requires this quantity to be small (see Appendix). In particular we implement \ref{synceq} as 
\[
(I-P_{n+1})\delta_{n+1} = (I-P_{n+1})[F(\bar{u}^{(k)}_n+(I-P_n)\delta_n) - \bar{u}^{(k)}_{n+1}].
\]
The initial condition $(\deltaStab)_0 = (I-P_0)\delta_0$ on time window $m$ is determined by imposing
\[
(I-P_0)\delta_0 = (I-P_0) [v_T - u_0^{(k)}]
\]
where $v_T$ is the converged iterate $u$ at the terminal time on the time window $m-1$. Effectively, by imposing continuity in the stable directions during the full assimilation, also across window boundaries, when solving \eqref{synceq} we can define a unique solution (this analysis point of view is also taken up in \cite{HaOlTi11} and \cite{ToHu13}).  To obtain a good initial condition for the algorithm we perform smoothing on an initialization window: i.e.~we employ the full Newton's algorithm (see \ref{sec:fullnewton}) on a short window and start the forced system \eqref{synceq} from there. This also improves the approximation of the unstable directions at the beginning of the window at which the projected method is started.

\section{Numerical experiments\label{sec:num}}

In the preceding sections we have outlined a data assimilation method based on a tangent space splitting into stable and non-stable subspaces.  As described, the method assumes noisy observations of the full state of the system (i.e.~observation operator $H$ the identity map on $\R^d$) at each time step, and no restrictions are placed on the length of the time interval.

In this section we demonstrate the behavior of the method for low dimensional test problems: the Lorenz models L63 and L96.  We study dependence on dimension of the projection operator and window lengths.  We compare the method with 4DVar, and investigate the approaches for incomplete observations and parameter estimation.  

In all experiments, the observations are generated from the truth by adding i.i.d.~zero-mean Guassian noise as in equation \ref{obs} with diagonal covariance matrix $E = \nu^2 I$, where $\nu^2$ denotes the variance of the noise process. As convergence criterion for the projected Newton's method we use that $\frac{\|\B{b}\|_2}{\|\B{u}\|_2}<10^{-15}$, where $\B{b}$ is the projected residual in \ref{reduced}.

\subsection{Dependence on projector in the L63 model}\label{Sec:L63}

The well-known Lorenz attractor \cite{Lo63} is a chaotic dynamical system commonly used as a test problem for data assimilation algorithms. The L63 model is
\begin{equation}\label{L63}
    \dot{x}_1 = \sigma(x_2-x_1), \quad 
    \dot{x}_2 = x_1(\rho-x_3)-x_2, \quad 
    \dot{x}_3 = x_1x_2-\beta x_3
\end{equation}
where $\sigma=10$, $\beta=\frac{8}{3}$ and $\rho=28$.  The Lyapunov exponents of the Lorenz attractor are $\lambda_1 \approx 0.906$, $\lambda_2 = 0$, $\lambda_3 \approx -14.572$.

For the experiments in this section we generate a (single) set of observations computing a trajectory of L63 on $t\in [0,20]$ with $T=20$, using time step $\Delta t = 0.005$, and $\nu^2 = 4$.
In all experiments in this section we use an assimilation window of length $\Delta \tau = 2.5$.

In \ref{sec:fullnewton} we observed that the full Newton's method successfully assimilates observations into L63. 
Now, we examine the proposed algorithm with projected Newton's and synchronization. 
Since the L63 model can be synchronized by coupling of the $x_1$-variables \cite{PeCa90,HaOlTi11}, it is natural 
instead of computing Lyapunov vectors to try to take $P=P_{x_1}$, hence always projecting on the $x_1$-coordinate, and to iterate \eqref{update1} and \eqref{synceq}. 
Errors~\eqref{MSE}--\eqref{anerror} are given in \ref{tabL63Px}, 
where it is clear that for our algorithm the choice $P=P_{x_1}$ is insufficient to obtain an orbit that is close to observations.   
Since the projection operators $P$ generally do not commute with the forward model solution operator $F$,
the projected Newton's method does not yield a projection of the full model solution, which means in particular that 
there are important differences between our algorithm and synchronization in the sense of \cite{PeCa90, HaOlTi11}.\par
Therefore we consider the projection operator on the subspace spanned by Lyapunov vectors. First, we choose the dimension of the projection operator to be $p=1$.
This means we use Newton's method in the (approximate) unstable direction, but not in the neutral or stable direction, because the L63 model has one positive, one zero and one negative Lyapunov exponent.
This is not sufficient for Newton's method to always converge, since the method works well until $t=20$ and after that Newton's method diverges. The results up to $t=20$ are shown in \ref{tabL63k1}, 
where in addition to errors we also display a measure of discontinuity at window boundaries defined as
\begin{equation}\label{jump}
 D = \frac{1}{N} \sum_{n=1}^N \max|G_n(\B{u})|
\end{equation}
and the average number of iterations needed for Newton's method to converge is denoted as $\#$.
We remark that in principle we could restart the method using full Newton's at $t=20$ and then continue with $p=1$.\par
\begin{figure}
\begin{center}
\includegraphics[width=0.8\textwidth]{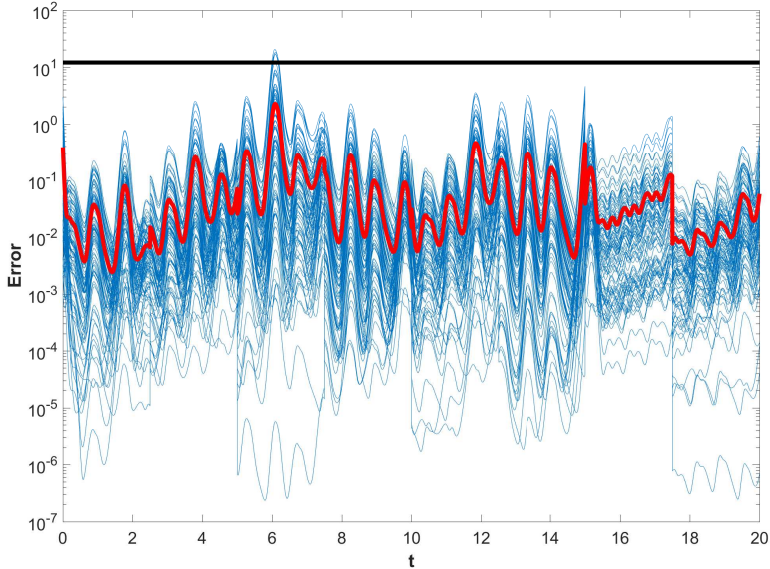}
\end{center}
\caption{Application of the algorithm to L63 with $P_2$. We compare the time-averaged observational error (black) with the error of the estimations (blue) over time. We used 100 observational noise realizations. The average estimation error over time is shown in red.}
\label{figL63k2}
\end{figure}

\begin{table}
\caption{Application of the algorithm to L63 with $P=P_{x_1}$. Results are unsatisfactory. Please recall $C$ and MSE are defined in equations \eqref{anerror} and \eqref{MSE} respectively.}
\label{tabL63Px}
\centering
\begin{tabular}{|l|c|}\hline
\textbf{Property}&\textbf{Value}\\\hline
Observation error $C(\left\lbrace X_n \right\rbrace)$&11.9\\
Distance between Estimate and Observations $C(\B{u})$&367\\
Error between Estimate and the Truth MSE&356\\\hline
\end{tabular}
\end{table}
\begin{table}
\caption{Application of the algorithm to L63 with $P_1$. Results for the time up to 20, since after that the algorithm diverges. Please recall $C$, MSE and $D$ are defined in equations \eqref{anerror}, \eqref{MSE} and \eqref{jump} respectively.}
\label{tabL63k1}
\centering
\begin{tabular}{|l|c|}\hline
\textbf{Property}&\textbf{Value}\\\hline
Observation error $C(\left\lbrace X_n \right\rbrace)$&11.9\\
Distance between Estimate and Observations $C(\B{u})$&12.2\\
Error between Estimate and the Truth MSE&0.27\\
Discontinuity measure at window boundaries D &0.23\\
Average number of iterations $\#$&$8.7$\\\hline
\end{tabular}
\end{table}
Next, we choose the dimension of the projection operator to be $p=2$. This means we apply Newton's method to both the unstable and neutral direction. 
The results are shown in \ref{tabL63k2,figL63k2}, where it can be seen that the algorithm becomes stable. Thus it is necessary to apply the projected Newton's method 
to both the unstable and neutral directions in this example. We repeated this numerical experiment for 100 different noise realizations.
\begin{table}
\caption{Application of the algorithm to L63 with $P_2$. Good results are obtained. Please recall $C$, MSE and $D$ are defined in equations \eqref{anerror}, \eqref{MSE} and \eqref{jump} respectively. Numbers shown in the table are averages over 100 noise realizations, together with the corresponding standard deviations.}
\label{tabL63k2}

\centering
\begin{tabular}{|l|c|}\hline
\textbf{Property}&\textbf{Value}\\\hline
Observation error $C(\left\lbrace X_n \right\rbrace)$&$12.00\pm 0.17$\\
Distance between Estimate and Observations $C(\B{u})$&$12.06\pm 0.18$\\
Error between Estimate and the Truth MSE&$0.09\pm 0.07$\\
Discontinuity measure at window boundaries D&$0.29\pm 0.08$\\
Average number of iterations $\#$&$6.52\pm 0.15$\\\hline
\end{tabular}
\end{table}
The case $p=d=3$ is the full Newton's method, which gives the smallest MSE as was shown in \ref{sec:fullnewton}.
We remark, however, that in this section initialization (full Newton's) is performed only on the first assimilation window, which reduces computational costs. 
We conclude that the algorithm is capable of recovering a good approximation of the true trajectory
and that this approximation \emph{is} a trajectory of the L63 model.\par
As mentioned in the beginning of \ref{sec:num}, observations of the full model state are not feasible and thus we need
to relax this assumption. Therefore we now assume that the only available observations are of the $x_1$-coordinate. 
First, we perform a preprocessing procedure in order to complete the missing observations: we run \eqref{receiver} with observations of $x_1$ as driving signal and $H \X_n=P_{x_1}\X_n$, $n=0,\dots,N$, as coupling. 
Subsequently, we apply the main algorithm with thus completed observations (which generally contain large errors due to the preprocessing).  
For the main algorithm we choose $p=2$. Results are shown in \ref{tabL63obsx}, where we see
that when only one coordinate is observed the error can be reduced and information on other coordinates can be obtained with synchronization as the preprocessing procedure.\par
\begin{table}
\caption{Application of the algorithm to L63 with $P_2$ and observations of the $x_1$-coordinate only. Synchronization is used as preprocessing step and errors are reduced by the projected Newton's method. Please recall $C$, MSE and $D$ are defined in equations \eqref{anerror}, \eqref{MSE} and \eqref{jump} respectively.}
\label{tabL63obsx}
\centering
\begin{tabular}{|l|c|}\hline
\textbf{Property}&\textbf{Value}\\\hline
Observation error $C(\left\lbrace X_n \right\rbrace)$&3.97\\
Distance between Estimate and Observations $C(\B{u})$&4.32\\
Error between Estimate and the Truth MSE&2.49\\
MSE for the observed coordinate&0.37\\
Discontinuity measure at window boundaries D&1.16\\
Average number of iterations $\#$&$7.0$\\\hline
\end{tabular}
\end{table}

\subsection{Dependence on window length in the L96 model}\label{Sec:L96}
Lorenz \cite{Lo96} proposed the following model as an example of a simple one-dimensional model with features of the atmosphere.  The L96 model is
\begin{equation}\label{L96}
\dot{x_l}=-x_{l-2}x_{l-1}+x_{l-1}x_{l+1}-x_l+\mathcal{F},\qquad (l=1,...,d),
\end{equation}
where the dimension $d$ and forcing $\mathcal{F}$ are parameters. Cyclic boundary conditions are imposed.
We implement the L96 model with the standard parameter choices $d=36$ and $\mathcal{F}=8$ . The differential equations are discretized with a forward Euler scheme with time step $\tau=0.005$
and the model initial conditions are chosen at random (standard Gaussian iid). 
Observations are obtained by perturbing a reference (true) trajectory with random Gaussian iid noise with zero mean and covariance $E=0.3^2I$. 
However, the observations are not drawn at every time step as for the L63 model but only every tenth time step, corresponding to observing a full model state every 6 hours. 
Then the map $F_n$~\eqref{map} corresponds to ten forward Euler steps.
This map is used to define $G$ and the derivatives of this map are needed for the QR-decompositions and Newton's iteration. 
For the synchronization we observe that if $G(\B{u})=0$, then we also have a trajectory under the forward Euler discretization with time step $\tau$. 
This means that any model integration can just be done with the original discretization, with forcing only applied at points where we have observations.\par
\begin{figure}
\begin{center}
\includegraphics[width=0.8\textwidth]{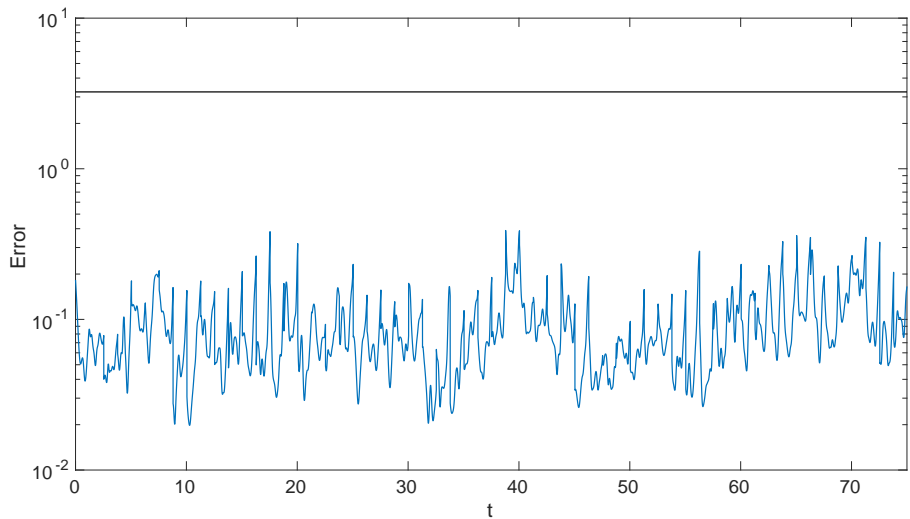}
\end{center}
\caption{The best estimation obtained by the algorithm applied to L96 ($p=15$, initialization window of 2.5, following windows of 1.25). }
\label{figL96best}
\end{figure}
\begin{figure}
\begin{center}
\includegraphics[width=0.8\textwidth]{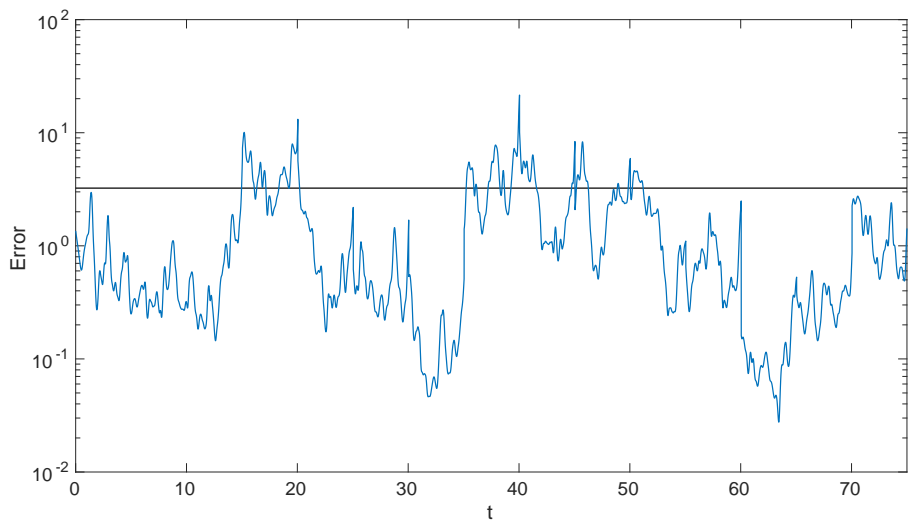}
\end{center}
\caption{The algorithm can also be applied to L96 with very long windows ($p=15$, initialization window of 15, following windows of 5). Trajectories are found, but if error reduction is the main purpose longer windows do not necessarily result in smaller errors.}
\label{figL96worst}
\end{figure}

For the projected Newton's method the dimension of the non-stable subspace $p$ is chosen to be either 15 or 25. 
We carry out numerical experiments for various choices of window lengths: we use initialization windows with lengths between 0.75 and 15 time units and following windows with lengths between 0.75 and 5 time units. 
The total time length of assimilation is always 75 and identical observations are used in all experiments.  
In \ref{figL96best,figL96worst} and \ref{L96tab15,L96tab25} 
it can be seen that the algorithm works well for both long and short windows, although when the windows are too long or too short the results deteriorate.
In general, higher $p$ decreases the estimation errors, although for the optimal choice of the window lengths---initialization window of 2.5 and following windows of 1.25---projection on $p=15$ results in better estimation, also see \ref{figL96best}.\par
\begin{table}
\caption{Application of the algorithm to L96 with $P_{15}$. Please recall $C$, MSE and $D$ are defined in equations \eqref{anerror}, \eqref{MSE} and \eqref{jump} respectively.}
\label{L96tab15}
\centering
\begin{tabular}{|l|c|c|c|c|c|c|c|c|}\hline
\textbf{Property}&\multicolumn{8}{|c|}{\textbf{Value}}\\\hline
\textbf{Window length}&\multicolumn{2}{|c|}{5}&\multicolumn{3}{|c|}{2.5}&\multicolumn{2}{|c|}{1.25}&0.75\\\hline
\textbf{Initial window}&15&5&15&5&\multicolumn{2}{|c|}{2.5}&1.25&0.75\\\hline
Observation error $C(\left\lbrace X_n \right\rbrace)$&\multicolumn{8}{|c|}{3.24}\\\cline{2-9}
Distance between Estimate and Observations $C(\B{u})$&4.70&3.65&3.58&3.36&3.35&3.22&3.25&3.26\\
Error between Estimate and the Truth MSE&1.47&0.41&0.40&0.18&0.18&0.09&0.12&0.20\\
Discontinuity measure at window boundaries D&0.96&0.51&0.48&0.35&0.35&0.21&0.24&0.28\\
Average number of iterations $\#$&11.1&10.4&8.8&8.5&8.6&7.5&7.5&7.0\\
\hline
\end{tabular}
\end{table}
\begin{table}
\caption{Application of the algorithm to L96 with $P_{25}$. Please recall $C$, MSE and $D$ are defined in equations \eqref{anerror}, \eqref{MSE} and \eqref{jump} respectively.}
\label{L96tab25}
\centering
\begin{tabular}{|l|c|c|c|c|c|c|c|c|}\hline
\textbf{Property}&\multicolumn{8}{|c|}{\textbf{Value}}\\\hline
\textbf{Window length}&\multicolumn{2}{|c|}{5}&\multicolumn{3}{|c|}{2.5}&\multicolumn{2}{|c|}{1.25}&0.75\\\hline
\textbf{Initial window}&15&5&15&5&\multicolumn{2}{|c|}{2.5}&1.25&0.75\\\hline
Observation error $C(\left\lbrace X_n \right\rbrace)$&\multicolumn{8}{|c|}{3.24}\\\cline{2-9}
Distance between Estimate and Observations $C(\B{u})$&4.01&3.94&3.38&3.26&3.26&3.16&3.16&3.09\\
Error between Estimate and the Truth MSE&0.82&0.74&0.22&0.11&0.12&0.09&0.10&0.16\\
Discontinuity measure at window boundaries D&0.89&0.87&0.37&0.33&0.32&0.27&0.27&0.30\\
Average number of iterations $\#$&9.1&9.1&8.0&8.0&8.0&7.0&7.0&6.9\\
\hline
\end{tabular}
\end{table}

In \ref{Lor95tab1} we investigate in more detail the dependence on the dimension of the projection $p$. We take total dimension $d=40$ and illustrate how the distance of the refined orbit from observations undergoes a sharp transition around the number of positive Lyapunov exponents, which is equal to 13. This sharp transition has also been observed for 4DVar-AUS \cite{TrDiTa10,PaCaTr13} The values that are reported are based upon uniformly distributed noise in the interval $[-2,2]$ and upon a total assimilation time of 4, subdivided into 4 windows of length 1. At the first window full Newton's is used (i.e. $p=40$) and at the subsequent three windows $p$ is as specified in \ref{Lor95tab1}; $C(\B{u})$ is then computed over all 4 windows. 

\begin{table}[!h]
\caption{Distance between Estimate and Observations $C(\B{u})$ for L96.}
\label{Lor95tab1}
\centering
\begin{tabular}{|l|c|c|c|c|c|c|c|}\hline
$p$ &  5 & 10 & 13 & 15 & 20 & 30 &40\\\hline
$C(\B{u})$ & 149.2 & 77.5 & 73.1 & 59.2 &  55.4 & 53.2 & 53.2\\\hline
\end{tabular}
\end{table}

It is possible to carry out the same experiment using more projected windows of length 1, to achieve a total length of 10 or 20. However, for small $p$, the results get worse if the total length increases. For length 50 divergence of Newton's method is observed. So, the table can be reconstructed qualitatively, but not for unlimited times. At some point new spin-up windows with full-Newton's ($p=40$) are needed. This eventual instability also occurs for L63 with $p=1$. This problem does not occur when $p$ is chosen large enough. Carrying out the numerical experiment with $p=20$ results in the average distance to observations remaining compatible with the noise level for time lengths up to 2500 (using windows of length 1).

This illustrates that it is important to define the non-stable space to be large enough, ensuring the $(I-P)$-problem does not contain neutral or unstable directions. If $p$ is chosen too small, the initialization at the spin-up window with full Newton's keeps the error somewhat in check over a few projected windows, but as we progress even further in time projected Newton's on an insufficiently large subspace is unable to keep the error in check. When errors get larger, this will eventually lead to divergence of Newton's method, but already before that the results from the data assimilation get progressively worse. However, if we choose p to be large enough, the method remains stable over long times.

\subsection{Comparison to 4DVar}\label{Sec:compto4dvar}
In the above sections we have argued that our algorithm aims at the same goals as the 4DVar algorithm and that for the L63 and L96 models we are able to 
reconstruct good trajectories based on observations. We now make a comparison with the standard 4DVar algorithm and demonstrate that our approach is a good alternative.\par
\begin{figure}
\begin{center}
\includegraphics[width=0.8\textwidth]{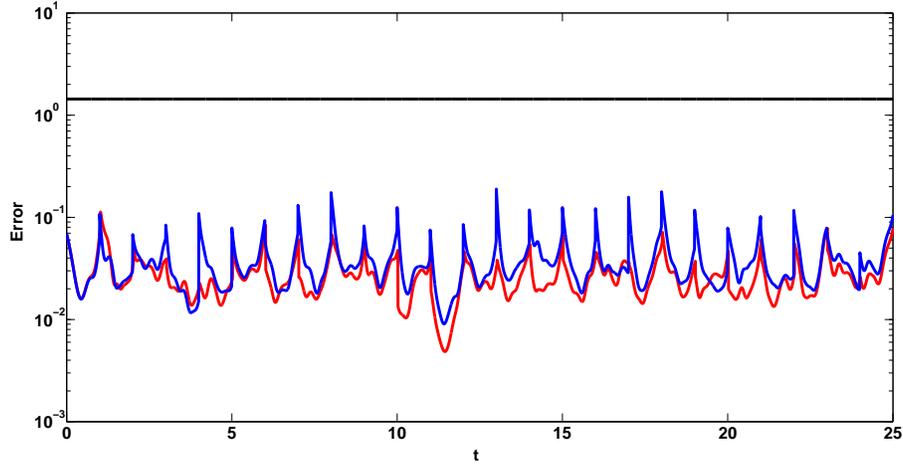}
\end{center}
\label{figerr4dvar}
\caption{The error between the results of the approximation methods and the truth over time, measured using the $\infty$-norm. The error of 4DVar is shown in blue, the error for our shadowing method is shown in green.}
\end{figure}

We perform a test using the L96 model with the same parameters as in the section above. Observations are drawn every fifth time step, which means we observe the full state 
every 0.025 time units, corresponding to 3 hours. We set $p=25$. In our tests we use identical data, models and windows for both methods. 
We choose 25 windows of length 1, of which the first is used as initialization window for the shadowing method. On the initialization window 5 iterations are needed for Newton's algorithm to converge. 
The initialization of 4DVar at the beginning of the first window is done with the first observation, since for neither of the two methods we have any prior knowledge of the system state. 
We do not use a background term for 4DVar. The gradient computation in
the 4DVar method is done using the adjoint integration and the
optimization is performed by a conjugate gradient
method. 
Some results are shown in \ref{figerr4dvar,tab4dvar}.

\begin{table}
\caption{The results from the projected Newton's algorithm and 4DVar are of comparable quality, but convergence is much more quick for the projected Newton's algorithm.}
\label{tab4dvar}
\centering
\begin{tabular}{|l|c|c|}\hline
\textbf{Property}&\multicolumn{2}{|c|}{\textbf{Value}}\\\hline
\textbf{Method}&Projected Newton's&4DVar\\\hline
Observation error $C(\left\lbrace X_n \right\rbrace)$&\multicolumn{2}{|c|}{1.43}\\\cline{2-3}
Distance between Estimate and Observations $C(\B{u})$ & 1.40 & 1.39\\
Error between Estimate and the Truth MSE&0.027&0.037\\
Discontinuity measure at window boundaries D&0.14&0.17\\
Average number of iterations $\#$&6.3&418.3\\\hline
\end{tabular}
\end{table}

\par

From this test it can be seen that 4DVar and our method give comparable results in this test case, but that our method is faster and more suitable if we were to use longer windows. The choice of window length is determined by the requirement that the 4DVar method should still converge. Even though the results in the previous section suggest slightly longer windows would be better, this is not viable for this test. In fact, the number of iterations could be significantly reduced for 4DVar by shortening the window, but this would come at the cost of not taking enough data into account.\par
We remark that one iteration of our method is more costly than one iteration of the 4DVar algorithm; if most directions are stable the difference in iteration cost would be less strong. In any case, the higher cost per iteration step of the projected Newton's method is more than compensated for since it requires far fewer iterations.  For sufficiently long windows, the cost per iteration of 4DVar is dominated by the need to do one model integration and one adjoint integration, which scales as $\mathcal{O}(Nc)$, where $c$ is defined as the typical cost of taking one time step in the non-linear model. An implementation of the full Newton's method in which the Jacobians are formed explicitly and are treated as dense matrices yields a cost per iteration of $\mathcal{O}(dNc)$ for the integrations needed and $\mathcal{O}(Nd^3)$ for solving the resulting linear system with a block tridiagonal method. The use of the projected Newton's method reduces this cost to $\mathcal{O}(pNc)+\mathcal{O}(Nd(p)^2)$. The main factors contributing to the cost per iteration of the projected Newton's method are two model integrations, $p$ tangent linear model integrations and the application of the modified Gram-Schmidt method to a $p\times d$-matrix at each time step.\par
We can see in \ref{tab4dvar} that the 4DVar method returns a result that is slightly closer to the observations than the projected Newton's method, while the projected Newton's method is slightly closer to the truth. This minor difference might be related to the projection on the unstable space used for the projected Newton's method \cite{TrDiTa10}.\par
Reasons for the large difference in number of iterations needed for convergence could be that Newton's method has a quadratic convergence rate, while the optimization algorithm for 4DVar does not. We remark that it is not possible to choose an algorithm with quadratic convergence rate for 4DVar, since we do not have the Hessian of the 4DVar cost function available. A more important reason for the strong difference in number of iterations needed for convergence may be in the fact that 4DVar and projected Newton's really solve very different problems. An explanation for the large difference in needed iterations for this example and the robustness of the projected Newton's method can be found by analyzing what happens when window lengths are increased.\par
The (projected) Newton's method has to solve larger (but weakly coupled and not that strongly nonlinear) problems if interval length is increased, while for the 4DVar approach the size of the optimization problem stays constant, but the problem becomes more and more highly nonlinear as the interval length is increased. This problematic behavior of the 4DVar optimization problem is well known in the literature \cite{Be91, MiGhGa94, PiVaTa96}. This can then lead to a large number of iterations needed for convergence, convergence to highly suboptimal local minima, or even to non-convergence of the 4DVar optimization. This difference between a large but weakly coupled and relatively easy root-finding problem that can be solved with an efficient method compared to a small but highly nontrivial optimization problem for which a slightly slower method has to be employed may give rise to the observed performance difference between the methods, both in terms of iterations needed to converge (and hence time needed to converge) and in the ability to still work for longer windows.

\subsection{Parameter estimation\label{sec:num_param_est}}
As described in \ref{sec:fullnewton}, shadowing-based data assimilation methodology can be applied to the problem of parameter estimation. 
The results of $\sigma$ estimation for the L63 model are shown in \ref{Globpartab},
where different values of initial $\sigma$ were chosen---5, 10, 15, and 20---with the true $\sigma$ being 10. 
Gaussian noise with identity covariance is added to the true solution and data assimilation is performed over one window of length 5 (when data assimilation is performed over multiple windows 
an estimate from the previous window can be taken as an initial parameter for the next window). 
It should be noted that similar results can be obtained for $\rho$ or $\beta$ estimation of the L63 model and for $\mathcal{F}$ estimation of the L96 model.\par
\begin{table}
\caption{Estimation of $\sigma$ by shadowing-based data assimilation methodology from \ref{sec:fullnewton}. The true value is 10.}
\label{Globpartab}
\centering
\begin{tabular}{|l|c|c|c|c|}\hline
\textbf{Property}&\multicolumn{4}{|c|}{\textbf{Initial guess for $\sigma$}}\\\hline
\textbf{Initial $\sigma$}&5&10&15&20\\\hline
Observation error $C(\left\lbrace X_n \right\rbrace)$&\multicolumn{4}{|c|}{2.97}\\\cline{2-5}
Distance between Estimate and Observations $C(\B{u})$ & 2.98 & 2.96 & 2.96 & 2.98\\
Error between Estimate and the Truth MSE& 0.03 & 0.02 & 0.03 & 0.07\\
Estimated $\sigma$ & 10.08 & 10.03 & 10.05 & 10.06\\\hline
\end{tabular}
\end{table}
In \ref{4DVarpartab2}, we show $\sigma$ estimations obtained by 4DVar using a window length of 0.25.
\begin{table}
\caption{Estimation of $\sigma$ by 4DVar parameter estimation. The true value is 10.}
\label{4DVarpartab2}
\centering
\begin{tabular}{|l|c|c|c|c|}\hline
\textbf{Property}&\multicolumn{4}{|c|}{\textbf{Initial guess for $\sigma$}}\\\hline
\textbf{Initial $\sigma$}&5&10&15&20\\\hline
Observation error $C(\left\lbrace X_n \right\rbrace)$&\multicolumn{4}{|c|}{2.97}\\\cline{2-5}
Distance between Estimate and Observations $C(\B{u})$ & 2.98 & 2.97 & 3.04 & 3.06\\
Error between Estimate and the Truth MSE& 0.12 & 0.12 & 0.17 & 0.18\\
Estimated $\sigma$ & 9.92 & 9.95 & 9.91 & 9.94\\\hline
\end{tabular}
\end{table}
Instead of using the method of \ref{sec:fullnewton}, we can also introduce trivial equations for the parameters to the shadowing-based data assimilation, which introduces  
extra zero Lyapunov exponents. As can be observed from \ref{Locparsigtab} this method fails for $\sigma$ estimations, 
though performs sufficiently well for $\rho$ (see \ref{Locparrhotab}) or $\beta$ estimations of the L63 model and for $\mathcal{F}$ estimation of the L96 model. 
Thus, adding trivial equations for the parameters to the shadowing-based data assimilation deteriorates its performance.
\begin{table}
\caption{Estimation of $\sigma$ by shadowing-based data assimilation methodology with trivial model for the parameters. 
The true value is 10 and the estimated $\sigma$ is the mean estimate. Cases when Newton's method diverges are denoted by ``$\infty$" in the corresponding column.}
\label{Locparsigtab}
\centering
\begin{tabular}{|l|c|c|c|c|}\hline
\textbf{Property}&\multicolumn{4}{|c|}{\textbf{Initial guess for $\sigma$}}\\ \hline
\textbf{Initial $\sigma$}&5&10&15&20\\\hline
Observation error $C(\left\lbrace X_n \right\rbrace)$ & \multicolumn{4}{|c|}{2.97} \\ \cline{2-5}
Distance between Estimate and Observations $C(\B{u})$ & 360 & 5.50 & $\infty$ & $\infty$\\
Error between Estimate and the Truth MSE& 355 & 2.77 & $\infty$ & $\infty$\\
Estimated $\sigma$ & 8.85 & 9.83 & $\infty$ & $\infty$\\ \hline
\end{tabular}
\end{table}
\begin{table}[!h]
\caption{Estimation of $\rho$ by shadowing-based data assimilation methodology with trivial model for the parameters. The true value is 28 and the estimated $\rho$ is the mean estimate.}
\label{Locparrhotab}\centering
\begin{tabular}{|l|c|c|c|c|} \hline
\textbf{Property}&\multicolumn{4}{|c|}{\textbf{Initial guess for $\rho$}}\\ \hline
\textbf{Initial $\rho$}&14&28&42&56\\\hline
Observation error $C(\left\lbrace X_n \right\rbrace)$&\multicolumn{4}{|c|}{2.97} \\ \cline{2-5}
Distance between Estimate and Observations $C(\B{u})$ & 35.7 & 2.95 & 5.52 & 161 \\
Error between Estimate and the Truth MSE& 32.7 & 0.02 & 2.74 & 159 \\
Estimated $\rho$ & 26.7 & 28.0 & 28.4 & 37\\ \hline
\end{tabular}
\end{table}

\section{Conclusions\label{sec:conclusions}}

We have introduced a new class of algorithms for data assimilation based upon shadowing refinement, synchronization, AUS, and PDA techniques.
Projections are determined based upon techniques employed in the computation of Lyapunov exponents/vectors, in particular continuous
QR techniques. This produces a splitting of the dynamics into non-stable and stable components,  
which allows for employing different techniques for the different components that are suited to their dynamics. Since the projections
are a function of solutions of the state space model, these projection based techniques require at least an approximate solution to 
determine initial projections. Assessing the uncertainty in obtaining an initial approximate solution and the impact of these uncertainties
on the assimilation is a focus of our future work. These techniques are also amenable in a number of ways to a Bayesian framework and since
we obtain an approximation of a time dependent orthonormal basis for the non-stable subspace one can assess the observation operator
with respect to the unstable subspace. The stable component has contractive dynamics which is useful for error control and further
assessing of uncertainties.
The algorithm developed here is effective in parameter estimation without introducing a trivial ODE for the parameters as in traditional data assimilation methods. 
We used a combination of analysis and numerical experiments to show that the algorithm works effectively and we demonstrated that the results compare favorably to those of 4DVar.
The other avenues for future research include more efficient numerical linear algebra techniques (the shadowing refinement relies on a block tridiagonal linear
system solve that we have performed with direct methods) and the use of parallel computing techniques.

\appendix

\section{Convergence of the synchronization update in the stable subspace\label{sec:analysis}}%
\subsection{Convergence in the linear, nonautonomous case}%
We study a synchronization process where there is some error made in the non-stable directions. If the model is linear but non-autonomous and at each step sufficiently close to the identity and the largest Lyapunov exponent of the stable subspace is negative, then the total error of the synchronized solution will not be much larger then the error in the non-stable directions. This holds in particular if the largest Lyapunov exponent of the stable subspace is small enough and if convergence to the Lyapunov exponents in the stable space is quick.\par
Let $\X_n$ be a solution to the nonautonomous linear model $\X_{n+1}=\D{F}_n\X_n$, for $n\in\mathbb{N}$, and 
let $Q_{n+1}R_{n+1}=\D{F}_nQ_{n}$.
Let $\bar{u}_n\in\R^d$ be a sequence of vectors approximating the truth in the non-stable subspace as follows
\[
	P_n \bar{u}_n=\bar{u}_n,  \qquad \|P_n \X_n-\bar{u}_n\|<\epsilon,
\]
where $P_n:=(Q_n^u)(Q_n^u)^T$, the orthogonal projector onto the first $p$ columns of $Q_n$, i.e.~$p$ is the dimension of the non-stable subspace.

Define $\Delta_n:=\left(\D{F}_n-I\right)$ and $\Delta:=\sup_n\|\Delta_n\|_F$.
Let $w_0$ be arbitrary and 
\[
	w_n=(I-P_n)\D{F}_{n-1}\left(\bar{u}_{n-1}+w_{n-1}\right),
\] 
for $n\ge 1$. Define the error vector $v_n$ as the difference between truth and approximation at time step $n$:
\begin{equation}\label{fullerr}
v_n:=\X_n-\bar{u}_n-w_n,
\end{equation}
and denote the projection on the $\kappa$-th column of $Q_n$ by $v_n^{\kappa}$, for $\kappa=p+1,\dots,d$.
Let $\zeta:=\left(\frac{\sqrt{2}\Delta}{1-\Delta}+\left[\frac{\sqrt{2}\Delta}{1-\Delta}\right]^2\right)(1+\Delta)\epsilon$ and $\delta^2=\left(2\frac{\sqrt{2}\Delta}{1-\Delta}\Delta+\left[\frac{\sqrt{2}\Delta}{1-\Delta}\right]^2(1+\Delta)\right)$. 
Define modified Lyapunov exponents as 
$\hat{\lambda}_k:=\lim_{n\to\infty}\frac{1}{n}\sum_{l=0}^{n}\log\left(R^{(k,k)}_\ell+\delta^2\right)$, for $k=p+1,\dots,d$.
\begin{assumption}\label{LEconv}
There exists a positive constant $\hat{\varepsilon}>0$ and a positive integer $\mathcal{N}$ such that for all integers $\mathcal{N}'\ge\mathcal{N}$ and all $n>0$,
\[
	\big| \frac{1}{\mathcal{N}'} \sum_{\ell=0}^{\mathcal{N}'} \log \left( R^{(k,k)}_{(n+\ell)} + \delta^2 \right) - \hat{\lambda}_k \big| < \hat\varepsilon, \quad k=p+1,\dots,d.
\]
\end{assumption}

\begin{theorem}\label{ThmLin}
Assume $\Delta<1$, $\hat{\lambda}_{p+1}<0$. Under \ref{LEconv}, for any $p<\kappa\leq d$, $m_1\in\mathbb{N}$ and $m_2:=m_1+\mathcal{N}+1$,
\begin{align*}
|v_{m_2}^{(\kappa)}|&\leq e^{(\hat{\lambda}_\kappa+\hat{\epsilon}_\kappa)\mathcal{N}}|v_{m_1}^{(\kappa)}|+
\sum_{l=0}^{\mathcal{N}}\zeta\prod_{j=l}^{\mathcal{N}-1}\left(R^{(\kappa,\kappa)}_{j+m_1+2}+\delta^2\right)\\
&\quad+\left\lbrace \sum_{\iota=\kappa+1}^d\sum_{l=0}^{\mathcal{N}}|R^{(\iota,\kappa)}_{l+m_1}+\delta^2|\prod_{j=l+1}^{\mathcal{N}}\left(R^{(\kappa,\kappa)}_{j+m_1}+\delta^2\right)|v_{l+m_1}^{(\iota)}|\right\rbrace.
\end{align*}
\end{theorem}
\begin{Corollary}\label{CorLin}
Assume that when averaging over all $n\in\mathbb{N}$, for all $m\in\mathbb{N}$ and all $\iota_l, \kappa_l \in \lbrace p+1,p+2,...,d\rbrace$, with $l\in\lbrace n,n+1,...,n+m\rbrace$, the average of the product $\overline{\Pi_{l=n}^{n+m}R^{(\iota_l,\kappa_l)}_l}$ can be expressed as the product of the averages: $\overline{\Pi_{l=n}^{n+m}R^{(\iota_l,\kappa_l)}_l}=\Pi_{l=n}^{n+m}\overline{R^{(\iota_l,\kappa_l)}_l}$.
Then 
\[
	\overline{|v_{m_2}^{(\kappa)}|}<\frac{1}{1-e^{(\hat{\lambda}_\kappa+\hat{\epsilon}_\kappa)\mathcal{N}}}\left(1+\frac{1}{|\hat{\lambda}_\kappa|}\left(1-e^{\hat{\lambda}_\kappa \mathcal{N}}\right)\right)\zeta+\mathcal{O}(\Delta^2),
\]
where $\overline{|v_{m_2}^{(\kappa)}|}$ denotes taking the average of $|v_{m_2}^{(\kappa)}|$ over all $m_2$.
\end{Corollary}
\noindent{\bf Proof of Theorem \ref{ThmLin}.}\par
To prove \ref{ThmLin}
we use, without loss of generality, a coordinate system such that at step $n$ the orthonormal Lyapunov vectors coincide with the standard basis.  
Then we consider the equation for the projection of the error on the stable space:
\begin{equation}\label{errorstepequation}
(I-P_{n+1})v_{n+1}=(I-P_{n+1})\left[Q_{n+1}R_{n+1}(\X_n-\bar{u}_n-w_n)\right].
\end{equation}
We now split the right side of \eqref{errorstepequation} between the range and the kernel of $P_n$:
\begin{align*}
(I-P_{n+1})&Q_{n+1}R_{n+1}(\X_n-\bar{u}_n-w_n) \\
=&(I-P_{n+1})Q_{n+1}R_{n+1}P_n(\X_n-\bar{u}_n-w_n) \\
&+(I-P_{n+1})Q_{n+1}R_{n+1}(I-P_n)(\X_n-\bar{u}_n-w_n)\\
=&(I-P_{n+1})Q_{n+1}R_{n+1}(P_n\X_n-\bar{u}_n)+(I-P_{n+1})Q_{n+1}R_{n+1}\left\lbrace(I-P_n)\X_n-w_n\right\rbrace.
\end{align*}
We first analyze the contribution from the error term $(I-P_{n+1})Q_{n+1}R_{n+1}(P_n\X_n-\bar{u}_n)$, which is the contribution of the error parallel to the $p$ leading Lyapunov vectors at time $n$ to the error perpendicular to these vectors at time $n+1$.\par
By the definitions of $Q_{n}$, $R_{n}$ and $\Delta_{n}$, we have that
\begin{equation}\label{QRpert}
Q_{n+1}R_{n+1}=I+\Delta_{n}.
\end{equation}
We recall that $P_n=Q_n^u(Q_n^ u)^T$ and that $Q_n$ is the identity matrix in our coordinates. We can approximate $||I-Q_{n+1}||_F\leq\frac{\sqrt{2}||\Delta_{n+1}||_F}{1-||\Delta_{n+1}||_2}$, by Theorem 3.1 of \cite{Ch12}. It immediately follows that
\begin{align*}
||(I-P_{n+1})Q_{n+1}R_{n+1}(P_n\X_n-\bar{u}_n)||_2&< \left(\frac{\sqrt{2}||\Delta_{n+1}||_F}{1-||\Delta_{n+1}||_2}+\left[\frac{\sqrt{2}||\Delta_{n+1}||_F}{1-||\Delta_{n+1}||_2}\right]^2\right)(1+||\Delta_{n+1}||_2)\epsilon\\
&<\zeta.
\end{align*}
\par
For the convergence in the stable directions we proceed analogously to \cite{VV01}. We remark that $R_n$ encodes the local approximation to the Lyapunov exponents \cite{ErPo98}. We recall that the Gram-Schmidt algorithm ensures that all diagonal elements of $R_n$ are positive for all $n$. Using the bound on the contribution of the unstable errors to the stable direction, we obtain
\begin{align}
||(I-P_{n+1})v_{n+1}||_2 &< \,||(I-P_{n+1})\D{F}_n\left\lbrace(I-P_n)v_n\right\rbrace||_2+\zeta\\\nonumber
\leq &\, ||(I-P_n)R_{n+1}\left\lbrace(I-P_n)v_n\right\rbrace||_2\\\nonumber
&+||((I-P_{n+1})Q_{n+1}-(I-P_n))R_{n+1}\left\lbrace(I-P_n)v_n\right\rbrace||_2+\zeta\\\nonumber
<&\, ||(I-P_n)R_n\left\lbrace(I-P_n)v_n\right\rbrace||_2\\\nonumber
+&\left(2\frac{\sqrt{2}||\Delta_{n+1}||_F}{1-||\Delta_{n+1}||_2}\Delta+\left[\frac{\sqrt{2}||\Delta_{n+1}||_F}{1-||\Delta_{n+1}||_2}\right]^2(1+\Delta)\right)||(I-P_n)v_n||_2+\zeta\\\nonumber
\leq & \, ||(I-P_n)R_n\left\lbrace(I-P_n)v_n\right\rbrace||_2+\delta^2||(I-P_n)v_n||_2+\zeta\nonumber
.
\end{align}
We now compute a bound in the expected error by induction on the stable subspace dimension $d-p$. If $d-p=1$,  then $(I-P_n)v_n=v_n^{(d)}$ and for any $m_1$, $m_2 \in \mathbb{N}$, with $m_2>m_1$,
\begin{equation}\label{LinSyncinitial}
|v_{m_2}^{(d)}|\leq\prod_{l=m_1}^{m_2-1}\left(\left(R_{dd}\right)_l+\delta^2\right)|v_{m_1}^{(d)}|+
\sum_{l=m_1+1}^{m_2}\zeta\prod_{j=l+1}^{m_2}\left(\left(R_{dd}\right)_j+\delta^2\right).
\end{equation}
Using \ref{LEconv} and choosing $m_2-m_1=\mathcal{N}+1$, we get
\begin{align*}
|v_{m_2}^{(d)}|&\leq\prod_{l=0}^{\mathcal{N}}\left(\left(R_{dd}\right)_{l+m_1}+\delta^2\right)|v_{m_1}^{(d)}|+
\sum_{l=0}^{\mathcal{N}}\zeta\prod_{j=l}^{\mathcal{N}-1}\left(\left(R_{dd}\right)_{j+m_1+2}+\delta^2\right)\\
&\leq e^{(\hat{\lambda}_d+\hat{\epsilon}_d)\mathcal{N}}|v_{m_1}^{(d)}|+
\sum_{l=0}^{\mathcal{N}}\zeta\prod_{j=l}^{\mathcal{N}-1}\left(\left(R_{dd}\right)_{j+m_1+2}+\delta^2\right)
\end{align*}
Now assume $d-p>1$ and let $p<\kappa\le d$
, then
\begin{align*}
|v_{m_2}^{(\kappa)}|&\leq\prod_{l=0}^{\mathcal{N}}\left(\left(R_{\kappa\kappa}\right)_{l+m_1}+\delta^2\right)|v_{m_1}^{(\kappa)}|+
\sum_{l=0}^{\mathcal{N}}\zeta\prod_{j=l}^{\mathcal{N}-1}\left(\left(R_{\kappa\kappa}\right)_{j+m_1+2}+\delta^2\right)\\
&\quad+\left\lbrace \sum_{\iota=\kappa+1}^d\sum_{l=0}^{\mathcal{N}}|(R_{\iota\kappa})_{l+m_1}+\delta^2|\prod_{j=l+1}^{\mathcal{N}}\left((R_{\kappa\kappa})_{j+m_1}+\delta^2\right)|v_{l+m_1}^{(\iota)}|\right\rbrace.
\end{align*}
We remark the first two terms of the above formula are the same as those in Eq.\eqref{LinSyncinitial}. This finishes the proof of the theorem. {\flushright $\square$\par}

Corollary \ref{CorLin} can be proven by taking the average over all $m_2$. 
We may now take averages over all $m_2>\mathcal{N}+1$ and use that $m_1=m_2-\mathcal{N}-1$.
\begin{align*}
\overline{|v_{m_2}^{(d)}|}&\leq e^{(\hat{\lambda}_d+\hat{\epsilon}_d)\mathcal{N}}\overline{|v_{m_2-\mathcal{N}-1}^{(d)}|}+
\sum_{l=0}^{\mathcal{N}}\zeta\prod_{j=l}^{\mathcal{N}-1}\left(\overline{\left(R_{dd}\right)_{j+m_2-\mathcal{N}+1}}+\delta^2\right)\\
&\leq e^{(\hat{\lambda}_d+\hat{\epsilon})\mathcal{N}}\overline{|v_{m_2}^{(d)}|}+
\sum_{l=0}^{\mathcal{N}}\zeta e^{\hat{\lambda}_d l}\\
&<e^{(\hat{\lambda}_d+\hat{\epsilon}_d)\mathcal{N}}\overline{|v_{m_2}^{(d)}|}+
\zeta\left(\int_{t=0}^{\mathcal{N}}e^{\hat{\lambda}_d t}\dd t+1\right)\\
&\leq e^{(\hat{\lambda}_d+\hat{\epsilon}_d)\mathcal{N}}\overline{|v_{m_2}^{(d)}|}+
\left(1+\frac{1}{|\hat{\lambda}_d|}\left(1-e^{\hat{\lambda}_d \mathcal{N}}\right)\right)\zeta.
\end{align*}
From which it immediately follows that
\begin{equation}
\overline{|v_{m_2}^{(d)}|}<\frac{1}{1-e^{(\hat{\lambda}_d+\hat{\epsilon}_d)\mathcal{N}}}\left(1+\frac{1}{|\hat{\lambda}_d|}\left(1-e^{\hat{\lambda}_d \mathcal{N}}\right)\right)\zeta.
\end{equation}
For the last term in the expression with $d-p>1$ this yields
\begin{align*}
&\qquad\overline{\sum_{\iota=\kappa+1}^d\sum_{l=0}^{\mathcal{N}}|(R_{\iota\kappa})_{l+m_1}+\delta^2|\prod_{j=l+1}^{\mathcal{N}}\left((R_{\kappa\kappa})_{j+m_1}+\delta^2\right)|v_{l+m_1}^{(\iota)}|}\\
&=\sum_{\iota=\kappa+1}^d\sum_{l=0}^{\mathcal{N}}\prod_{j=l+1}^{\mathcal{N}}\left(\overline{(R_{\kappa\kappa})_{j+m_1}}+\delta^2\right)\overline{|(R_{\iota\kappa})_{l+m_1}+\delta^2|}\overline{|v_{l+m_1}^{(\iota)}|}\\
&<\sum_{\iota=\kappa+1}^d \left(1+\frac{1}{|\hat{\lambda}_\kappa|}\left(1-e^{\hat{\lambda}_\kappa \mathcal{N}}\right)\right)\left(\Delta+\delta^2\right)\overline{|v_{m_2}^{(\iota)}|}.
\end{align*}
From which it immediately follows that
\begin{equation}\label{LinSyncinduction}
\overline{|v_{m_2}^{(\kappa)}|}<\frac{1}{1-e^{(\hat{\lambda}_\kappa+\hat{\epsilon}_\kappa)\mathcal{N}}}\left(1+\frac{1}{|\hat{\lambda}_\kappa|}\left(1-e^{\hat{\lambda}_\kappa \mathcal{N}}\right)\right)\left(\zeta+\sum_{\iota=\kappa+1}^d\Delta\overline{|v_{m_2}^{(\iota)}|}\right).
\end{equation}
Combining Eqs. \eqref{LinSyncinitial} and \eqref{LinSyncinduction} and using that $\zeta=\mathcal{O}(\Delta)$, we conclude that
\begin{equation}
\overline{|v_{m_2}^{(\kappa)}|}<\frac{1}{1-e^{(\hat{\lambda}_\kappa+\hat{\epsilon}_\kappa)\mathcal{N}}}\left(1+\frac{1}{|\hat{\lambda}_\kappa|}\left(1-e^{\hat{\lambda}_\kappa \mathcal{N}}\right)\right)\zeta+\mathcal{O}(\Delta^2).
\end{equation}
{\flushright $\square$\par}

\subsection{Bound for the nonlinear case}
For the nonlinear case we do not have a convergence proof, but we can put a bound on the error. Let the truth $\X_n$ be a solution to the nonlinear model $\X_{n+1}=F_n(\X_n)$, where $F_n$ is a $\mathcal{C}^3$ function. Assume $\{ \X_n \}$ lies on an attractor of $F$ and on the attractor $F$ admits an exponential splitting. Let $\epsilon_1\geq\epsilon_2>0$, $A_{\epsilon_1}$ the neighborhood of size $\epsilon_1$ around the attractor of $F$, $\alpha\geq 0$, $\delta>0$ and $\tilde{\lambda}>\exp(\lambda_s)$, where $\lambda_s<0$ is the largest Lyapunov exponent of the stable space. Let $\Pi_n$ be projectors that project on the non-stable space at $\X_n$, let $K_2=\frac{1}{2}\sup_{\chi\in A_{\epsilon_1} } |\DD{F}_n(\chi)|$ and $K_0=\sup_{\chi\in A_{\epsilon_1} } |F_n(\chi)|$. Let $P_n\in\R^d\times\R^d$ be a sequence projectors and let $\bar{u}_n\in\R^d$ be a given sequence of vectors with $P_n \bar{u}_n = \bar{u}_n$ for all $n$. Let $w_0$ be some arbitrary vector and $w_n=(I-P_n)F\left(\bar{u}_{n-1}+w_{n-1}\right)$, for $n>1$. Define the error vector $v_n$ as the difference between truth and approximation at time step $n$:
\begin{equation*}
v_n:=\X_n-\bar{u}_n-w_n.
\end{equation*}
\begin{theorem}\label{ThmNonlin}
Assume $||v_n||<\epsilon_1$, $||P_{n}v_{n}||<\epsilon_2$ and $||P_{n+1}v_{n+1}||=||P_{n+1}\X_{n+1}-\bar{u}_{n+1}||<\epsilon_2$. Then there exists some $\tilde{\alpha}>0$ such that if for all $\tilde{v}\in\R^d$ it holds that $P_n \tilde{v}\in K^u_\alpha(\X_n)$, where $K^u_\alpha(\X_n)$ is the non stable cone \cite{KaHa95, BrSt02} of size $\alpha$ at $\X_n$, and $(I-P_n)\tilde{v}\in K^s_{\tilde{\alpha}}(\X_n)$, where $K^s_{\tilde{\alpha}}(\X_n)$ is the stable cone \cite{KaHa95, BrSt02} of size $\tilde{\alpha}$ at $\X_n$, then
\begin{equation}
||v_{n+1}||<\epsilon_2+K_2\epsilon_1^2+2\alpha (K_0\epsilon_2+K_2\epsilon_2^2)+(\tilde{\lambda}+\delta)\epsilon_1.
\end{equation}
\end{theorem}
\noindent{\bf Proof of Theorem\ref{ThmNonlin}.}\par
Throughout this proof we use results of \cite{KaHa95, BrSt02, ToHu13}.
Since the error $v_n$ is small and $F_n$ is $\mathcal{C}^3$, we can approximate the nonlinear flow by a Taylor expansion around the truth
\begin{equation}\label{texp}
||v_{n+1}-\D{F}_n(\X_n)v_n||=\|F_n(\X_n)-F_n(\X_n-v_n)-\D{F}_n(\X_n) v_n \|\leq K_2||v_n||^2.
\end{equation}
By splitting $v_{n+1}=P_{n+1}v_{n+1}+(I-P_{n+1})v_{n+1}$, noting that $||I-P_{n+1}||<1$ and using \eqref{texp} we obtain 
\begin{equation}\label{spliterr}
||v_{n+1}||<\epsilon_2+K_2\epsilon_1^2+||(I-P_{n+1})\D{F}_n(\X_n)v_n||.
\end{equation}
Due to the exponential splitting of $F_n$, the non-stable cone becomes more narrow under the tangent dynamics, i.e. vectors in 
this non-stable cone tend to align more towards the non-stable directions under the dynamics. 
This means that $\D{F}_n(\X_n)K^u_\alpha(\X_n)\subset\mathrm{int}(K^u_\alpha(\X_{n+1}))\cup \{0\}$ (it follows for example from Proposition 5.4.1 of \cite{BrSt02} or Lemma 6.2.10 of \cite{KaHa95}). 
Hence we have that $\D{F}(\X_n)P_nv_n\in\mathrm{int}(K^u_\alpha(\X_{n+1}))\cup \{0\}$. Due to the dynamics on the non-stable cone, we expect the length of $P_nv_n$ to grow.
A bound for the growth in any step is given by Taylor expansion of $F_n(\X_n-P_nv_n)$ around $F_n(\X_n)$ as $||\D{F}(\X_n)P_nv_n||\leq K_0||P_nv_n||+K_2||P_nv_n||^2$, 
where $K_0=\sup |F_n(\X_n)|$. However, the only part of $\D{F}(\X_n)P_nv_n$ of interest is the component $(I-P_{n+1})\D{F}(\X_n)P_nv_n$. 
We have that 
\begin{equation}\label{projerr}
||(I-P_{n+1})\D{F}(\X_n)P_nv_n||< 2\alpha ||\D{F}(\X_n)P_nv_n||<2\alpha (K_0||P_nv_n||+K_2||P_nv_n||^2).
\end{equation}
For the part $(I-P_n)v_n$ in the stable cone we can use \cite{BrSt02} Proposition 5.4.2 or \cite{KaHa95} Lemma 6.2.11, which states that this vector shrinks under time evolution, 
where the amount depends on the width of our cone. To be precise: $\forall\delta>0\;\exists\tilde{\alpha}>0$ such that if 
$(I-P_n)v_n\in K^s_{\tilde{\alpha}}(\X_n)$, then
\begin{equation}\label{flowerr}
\|\D{F}_n(\X_n)(I-P_n)v_n\|<(\tilde{\lambda}+\delta)\|(I-P_n)v_n\|. 
\end{equation}

Collecting the estimates~\eqref{spliterr}--\eqref{flowerr}, we find that
\begin{align*}
||v_{n+1}||&<\epsilon_2+K_2\epsilon_1^2+||(I-P_{n+1})\D{F}_n(\X_n)(P_nv_n+(I-P_n)v_n)||\\
&\leq \epsilon_2+K_2\epsilon_1^2+||(I-P_{n+1})\D{F}_n(\X_n)P_nv_n||+||\D{F}_n(\X_n)(I-P_n)v_n)||\\
&<\epsilon_2+K_2\epsilon_1^2+2\alpha (K_0||P_nv_n||+K_2||P_nv_n||^2)+(\tilde{\lambda}+\delta)||(I-P_n)v_n||\\
&<\epsilon_2+K_2\epsilon_1^2+2\alpha (K_0\epsilon_2+K_2\epsilon_2^2)+(\tilde{\lambda}+\delta)\epsilon_1.
\end{align*}
{\flushright $\square$\par}

\section*{Acknowledgments}
Andrew Steyer and Erik Van Vleck acknowledge support from NSF grant DMS-1419047. Xuemin Tu acknowledges support from NSF grant DMS-1419069. Bart de Leeuw acknowledges this work is part of the research programme Mathematics of Planet Earth 2014 EW with project number 657.014.001, which is financed by the Netherlands Organisation for Scientific Research (NWO).
\bibliographystyle{plain}
\bibliography{artikel}

\end{document}